\documentclass[final,leqno,onefignum,letterpaper]{siamltex704}

\usepackage{epsfig,amsfonts,latexsym,amsmath,graphicx,color}

\usepackage[ruled]{algorithm2e}

\usepackage{stackrel}

\usepackage{soul}

\DeclareGraphicsExtensions{.pdf,.eps}

\definecolor{bsb}{rgb}{0.1,0.1,0.7}

\definecolor{ggg}{rgb}{0.1,0.5,0.1}
\definecolor{rrr}{rgb}{0.8,0.2,0.2}

\newcommand{\pdpd}[2]{\frac{\partial #1}{\partial #2}}
\newcommand{\pdpdo}[2]{\frac{{\mathrm{d}} #1}{{\mathrm{d}} #2}}

\newcommand{\dd}{\mathrm{d}}
\newcommand{\ii}{\mathrm{i}}
\newcommand{\R}{\mathbb{R}}

\newcommand{\half}{\mbox{$\frac{1}{2}$}}

\newcommand{\ignore}[1]{}

\newcommand{\centerPar}[1]{\xi'_{#1}}

\newcommand{\lengthPar}[1]{L'_{#1}}
\newcommand{\lengthPer}[1]{L''_{#1}}
\renewcommand{\ii}{\mathrm{i}}
\newcommand{\ZZ}{{\mathbb Z}}

\newcommand{\Tvk}{T_{\nu,k}}

\SetAlFnt{\small} \SetAlCapFnt{\small} \SetAlCapNameFnt{\small}
\SetAlCapHSkip{0pt} \IncMargin{-\parindent}

\newcommand{\notmid}{\mid\kern-0.5em\not\kern0.5em}

\newcommand{\mstrut}[1]{\mbox{\rule{0mm}{#1}}}

\newcommand{\ba}{\begin{array}}
\newcommand{\ea}{\end{array}}

\title{Multiscale reverse-time-migration-type imaging \\ using
       the dyadic parabolic decomposition \\ of phase space
               }

\author{Fredrik Andersson\thanks{Mathematics LTH, Centre for
        Mathematical Sciences, Lund Institute of Technology, Lund
        University, Sweden ({\tt fa@maths.lth.se})}
\and
 Maarten V. de Hoop\thanks{Department of Mathematics, Purdue
   University, West Lafayette, IN 47907, USA ({\tt
     mdehoop@purdue.edu}).}
\and
 Herwig Wendt\thanks{CNRS, Institut de Recherche en Informatique de Toulouse (IRIT), Toulouse, France ({\tt herwig.wendt@irit.fr})}
        }

\begin{document}


\maketitle

\setcounter{page}{1}

\begin{abstract}
We develop a representation of reverse-time migration in terms of Fourier integral operators the  canonical relations of which are graphs. Through the dyadic parabolic decomposition of phase space, we obtain the solution of the wave equation with a boundary source and homogeneous initial conditions using wave packets. On this basis, we develop a numerical procedure for the reverse-time continuation from the boundary of scattering data and for RTM migration. The algorithms are derived from those we recently developed for the discrete approximate evaluation of the action of Fourier integral operators and inherit from their conceptual and numerical properties.
\end{abstract}

\pagestyle{myheadings}
\thispagestyle{plain}
\markboth{ANDERSSON ET AL.}{MULTISCALE RTM-TYPE IMAGING}

\section{Introduction}

Reflection seismology is a commonly used method to study the
properties of Earth's subsurface in geophysical exploration. Point
sources are placed on Earth's surface which generate acoustic waves in
the subsurface that are reflected where the medium properties vary
discontinuously. These reflections are recorded at Earth's surface by
(arrays of) point receivers. The goal of seismic imaging is to
reconstruct the singular variations in medium properties from the
reflected waves recorded at the surface
\cite{Claerbout1985Book,Biondi2006,Symes2009}. The most common
formulation for seismic inverse scattering takes the form of a
linearized inverse problem for the medium coefficient in the acoustic
wave equation, where the linearization is performed about a smoothly
varying background. Here, the background model is assumed to be
known. However, via a formulation as a separable inverse problem one
can also proceed with determining this background model. The
linearization defines a single scattering operator that maps the
coefficient contrast to the data, i.e., the restriction of the
scattered wave field to the acquisition set. The adjoint defines the
process of imaging.

We consider reverse-time migration (RTM) \cite{SchultzSherwood1980,
  whitmore:1983, mcmechan:1983, baysalKS:1983, sunM:2001}, and the
RTM-based inverse scattering transform developed and analyzed by Op't
Root \textit{et al.} \cite{OptrootSdH:2012}. Through an appropriate
formulation based on wave field continuation
\cite{Duchkov-wavemotion}, we arrive at a representation of RTM in
terms of a Fourier Integral Operator (FIO) associated with a canonical
graph. Indeed, we use such a representation. The key contribution of
this paper lies in the development of an algorithm for the
solution of the wave equation with a \textit{boundary source} and
homogeneous initial conditions using the dyadic parabolic
decomposition of phase space in the limit of fine scales. This
algorithm is then composed with an imaging condition to yield the
RTM-type imaging or inverse scattering. We explicitly admit the formation of caustics.

Viewing wave packets as
localized plane waves, our approach has connections to methods
in which one designs sources that favor (directional) illumination of
particular subdomains of the subsurface. We mention plane-wave
migration and beam-wave migration. In plane-wave migration one
synthesizes plane-wave source experiments \cite{Whitmore1995}. Given a plane-wave source one can then
introduce tilted coordinates to carry out the wave field extrapolation with a limited accuracy propagator \cite{ShanBiondi2004}. In beam-wave migration, Brandsberg-Dahl
and Etgen \cite{BrandsbergEtgen2003} use a rotating coordinate system and
essentially couple wave field methods with band-limited properties to
ray-geometric methods. Furthermore, we mention the use of coherent
states in this context by Albertin \textit{et al.} \cite{AlbertinYingst2002}. Instead of tilted
coordinates, one can use curvilinear coordinates in combination with a
paraxial propagator \cite{SavaFomel2005}; the curvilinear
coordinates may be generated as geodesics initiated from a point
source or a plane wave. We note, however, that these methods are
downward-continuation based whereas our approach is based on
reverse-time continuation; also, we decompose the data although we
could incorporate a synthesis of wave packet sources as well.

Our numerical solution is derived from the algorithm that we developed
for FIOs \cite{AndersondHW:2012}. Computational
efficiency arises from organizing the decomposition and propagation by
directions associated with frequency boxes instead of individual wave
packets. The superposition of wave packets is complete and their
propagation, as well as the corresponding imaging, converges in the
limit of fine scales in smooth velocity models.
Our formulation is insensitive to specific choices of (absorbing)
boundary conditions, which is in contrast to PDE based solutions,
including beam-wave migration. Moreover, it naturally conveys angular
information which can be used in the imaging process, for instance,
for computing restricted angle transforms.
Cand\`{e}s, Demanet \& Ying developed a fast butterfly algorithm for
the application of Fourier integral operators associated with
canonical graphs \cite{Candes2009}, which presents an interesting
alternative to the propagation component of our algorithm.
%

Our algorithm is particularly well suited for application to (limited
aperture) seismic array data providing a way of partial (in phase
space) imaging possibly with a small set of sources. Moreover, if one
need not generate a `global' image, we do not need to evaluate the
relevant wave field solutions at `all' times, unlike algorithms based
on numerically solving the wave equation, which enables
computationally efficient target-oriented imaging. Target-oriented
imaging can be used effectively, for example, with available arrays
and earthquakes in studying heterogeneities and discontinuities in
Earth's mantle \cite{Li2013, Vinnik2000}.
A key element of our algorithm is finding a low-rank separated
representation of the amplitude of the relevant FIO, which we do using Prolate Spheroidal Wave Functions (PSWFs) \cite{AndersondHW:2012}. Demanet \&
Ying \cite{Demanet2012} proposed a method of finding such a
representation based on the randomized sampling algorithm for
constructing factorizations for low-rank matrices.

Our algorithm involves the propagation of high-frequency waves. To
compare its complexity with the computational complexities of RTM
algorithms based on numerically solving the wave equation, one can
essentially compare the complexities of the backpropagation of a
boundary source (from single-source data). 
Considering (backward) solving the wave equation in dimension $n$ on a grid of side length $N$ and a propagation time of order $\mathcal{O}(1)$, the CFL condition implies that the time step is of order $\mathcal{O}(1/N)$. The $\mathcal{O}(N)$ applications of a short-time propagator, with a presumed complexity $\mathcal{O}(N^n)$, then yields a complexity $\mathcal{O}(N^{n+1})$. 
The time step in our algorithm is $\mathcal{O}(1)$ in principle, while the application of the Fourier integral operator representing the propagator is $\mathcal{O}(N^n \log N)$ per frequency box (see \cite{AndersondHW:2012}); the number of time steps to be computed is
affected by the size of the target. 
Demanet \& Ying \cite{Demanet2012} already pointed out the time upscaling of an
approach to wave propagation using a FIO.

The outline of the paper is as follows. In Section~\ref{sec:2}, we
summarize the parametrix construction of the wave equation, and
introduce the relevant Hamilton system and linearized Hamilton-Jacobi
equations describing the geometry of the imaging process. In
Section~\ref{sec:3} we formulate reverse-time continuation from the
boundary and obtain a particular oscillatory integral representation
for the kernel of this process, to which the algorithm for FIOs that we developed in an earlier paper applies.
Subsection~\ref{sec:33} contains this key new result, and
Subsection~\ref{sec:comprtc} its computational counterpart. These are
also the main components of the asymptotic form of the RTM-based
inverse scattering transform and imaging algorithm, which we develop
in Section~\ref{sec:4}.
In Section~\ref{sec:5}, we give
numerical examples of reverse-time continuation and inverse scattering
also in the presence of caustics. We end with a discussion in
Section~\ref{sec:6}.

\subsection*{Dyadic parabolic decomposition of phase space}

We briefly discuss the (co)frame of curvelets and wave packets
\cite{Cand`es2006,Duchkov2009a,Smith1998}. We will implicitly suppose that the data are decomposed into wave packets below, and we will develop wave packet based algorithms with accuracy $\mathcal{O}(2^{-k/2})$ \cite{AndersondHW:2012}.
\\
Let $u \in L^2(\R^n)$ represent a (seismic) velocity field, and let $\hat{u}(\xi) = \int u(x) \, \exp[-\ii \langle x,\xi  \rangle ] \, \dd x$ be  the Fourier
transform.
One begins with covering the positive $\xi_1$-axis
($\xi' = \xi_1$) by overlapping boxes of the form
\begin{equation}
\label{equ:box}
   B_k
   = \left[\centerPar{k} - \frac{\lengthPar{k}}{2},
           \centerPar{k} + \frac{\lengthPar{k}}{2}\right] \times
     \left[-\frac{\lengthPer{k}}{2},
            \frac{\lengthPer{k}}{2}\right]^{n-1}.
\end{equation}
Here, both the centers $\centerPar{k}$ and the side lengths
$\lengthPar{k}$, $\lengthPer{k}$ follow parabolic scaling
\[
   \centerPar{k} \sim 2^k,\quad \lengthPar{k} \sim 2^k ,
   \quad \lengthPer{k} \sim 2^{k/2} ,\quad
                        \mbox{as $k \rightarrow \infty$} .
\]
Next, for each $k \ge 1$, let $\nu$ vary over a set of approximately
$2^{k (n-1)/2}$ uniformly distributed unit vectors\footnote{By
convention, we let $\nu(0) = e_1$ be aligned with the $\xi_1$-axis.}. Let
$\Theta_{\nu,k}$ denote a choice of rotation matrix which maps $\nu$
to $e_1$, and
\begin{equation}
\label{equ:rotbox}
   B_{\nu,k} = \Theta_{\nu,k}^{-1} B_k .
\end{equation}
In the (co)frame construction, one
encounters two sequences of smooth functions, $\hat\chi_{\nu,k}$ and
$\hat\beta_{\nu,k}$, on $\R^n$, each supported in $B_{\nu,k}$, that form a copartition of unity
\begin{equation} \label{eq:3.1}
   \hat{\chi}_0(\xi) \hat{\beta}_0(\xi) + \sum_{k \ge 1} \sum_{\nu}
   \hat{\chi}_{\nu,k}(\xi) \hat{\beta}_{\nu,k}(\xi) = 1 ,
\end{equation}
and satisfy the estimates
\[
   |\langle \nu,\partial_{\xi} \rangle^j \, \partial_{\xi}^{\alpha}
    \hat{\chi}_{\nu,k}(\xi)| +
   |\langle \nu,\partial_{\xi} \rangle^j \, \partial_{\xi}^{\alpha}
    \hat{\beta}_{\nu,k}(\xi)|
   \le C_{j,\alpha} \, 2^{-k(j + |\alpha|/2)} .
\]
One now forms
\begin{equation} \label{eq:3.2}
   \hat{\psi}_{\nu,k}(\xi) = \rho_k^{-1/2}
             \hat{\beta}_{\nu,k}(\xi)\,,\quad
   \hat{\varphi}_{\nu,k}(\xi) = \rho_k^{-1/2}
              \hat{\chi}_{\nu,k}(\xi) ,
\end{equation}
where $\rho_k$ is the volume of $B_k$. These functions satisfy the
estimates
\begin{equation}
\label{eq:cdest}
  \forall N:\quad \left.\begin{array}{l} |\varphi_{\nu,k}(x)| \\[0.3cm]
                |\psi_{\nu,k}(x)| \end{array} \right\}
          \le C_N 2^{k (n+1)/4} \, (\,
       2^k |\langle\nu,x\rangle| + 2^{k/2} \| x \| \,)^{-N}.
\end{equation}
To obtain a (co)frame, one introduces the integer
lattice: $X_j := (j_1,\dots,j_n)\in\ZZ^n$, the dilation matrix
\[
   D_k = \frac{1}{2 \pi}
         \left(\begin{array}{lr} \lengthPar{k} & 0_{1 \times n-1} \\
     0_{n-1 \times 1} & \lengthPer{k} I_{n-1} \end{array}\right) ,
   \qquad \det \, D_k = (2\pi)^{-n} \rho_k ,
\]
and the points $x^{\nu,k}_j = \Theta_{\nu,k}^{-1} D_k^{-1} X_j$. The frame elements are now defined
in the Fourier domain as
\begin{equation} \label{eq:3.3}
   \hat{\varphi}_{\gamma}(\xi)
   = \hat{\varphi}_{\nu,k}(\xi) \,
         \exp[-\ii \langle x^{\nu,k}_j,\xi \rangle] ,\quad
     \gamma = (j,\nu,k) ,\quad k \ge 1,
\end{equation}
and similarly for $\hat{\psi}_{\gamma}(\xi)$. Thus, one obtains the transform pair
\begin{equation} \label{eq:ctrp}
   u_{\gamma} = \int u(x) \overline{\psi_{\gamma}(x)} \, \dd x ,
\quad\quad
   u(x) = \sum_{\gamma} u_{\gamma} \varphi_{\gamma}(x)
\end{equation}
with the property that $\sum_{\gamma':\ k'=k,\ \nu'=\nu} \!
u_{\gamma'} \hat{\varphi}_{\gamma'}(\xi) = \hat{u}(\xi)
\hat{\beta}_{\nu,k}(\xi) \hat{\chi}_{\nu,k}(\xi)$ for each $\nu, k$.

\section{Parametrix}
\label{sec:2}

Here, we summarize the parametrix construction for the wave equation. We
consider the Cauchy initial value problem
\begin{eqnarray}
   \left[\pdpd{{}^2}{t^2} + A(x,D_x)\right] u &=& 0 ,\quad
   A(x,D_x) = c(x) D_x^2 c(x) ,
\\
   u(x,0) &=& 0 ,\quad \pdpd{u}{t}(x,0) = h(x) ;
\end{eqnarray}
we have normalized the pressure: $u = c^{-1} p$.

To evaluate the parametrix, we use the first-order system for $u$ that
is equivalent to this wave equation,
\begin{equation}
\label{eq:first_order_time}
   \pdpd{}{t} \left( \ba{cc} u \\ \pdpd{u}{t}
              \ea \right)
   = \left( \ba{cc} 0 & 1 \\ - A(x,D_x) & 0 \ea \right)
    \left( \ba{cc} u \\ \pdpd{u}{t} \ea \right) .
\end{equation}
This system can be decoupled, namely, by the matrix-valued
pseudodifferential operators
\begin{equation*}
   V(x,D_x) =
   \left( \ba{cc} 1 &  1 \\
             -\ii B(x,D_x)\ & \ii B(x,D_x)\ \ea \right),\;\;
   \Lambda(x,D_x)
   = \half \left( \ba{cr} 1 & \ii B(x,D_x)^{-1} \\
             1 & -\ii B(x,D_x)^{-1} \ea \right) ,
\end{equation*}
where $B(x,D_x) = \sqrt{A(x,D_x)}$ is a pseudo\-differential operator
of order 1. 

The use of a general symbol $B$ in our presentation facilitates the
extension of our algorithm to the imaging with elastic
waves \cite{Brytik2013}.

The principal symbol of $B(x,D_x)$ is given by $B^{\rm
  prin}(x,\xi) = \sqrt{A^{\rm prin}(x,\xi)}$. Then
\begin{equation}
\label{eq:decoupled_firstorder_u_f}
   u_{\pm} = \half u \pm
         \half \ii B(x,D_x)^{-1} \pdpd{u}{t} ,
\end{equation}
satisfy the two first-order (``half wave'') equations
\begin{equation}
\label{eq:decoupled_firstorder}
   P_{\pm}(x,D_x,D_t) \, u_{\pm} = 0 ,
\end{equation}
where 
\begin{equation}
   P_{\pm}(x,D_x,D_t) = \pdpd{}{t} \pm \ii B(x,D_x) ,
\end{equation}
supplemented with the initial conditions
\begin{equation}
\label{eq:hpolIVs}
   u_{\pm} |_{t=0} = h_{\pm} ,\quad
   h_{\pm} = \pm \half \ii B(x,D_x)^{-1} h .
\end{equation}

We construct operators $S_{\pm}(t)$ that solve the initial value
problem (\ref{eq:decoupled_firstorder}), (\ref{eq:hpolIVs}):
$u_{\pm}(y,t) = (S_{\pm}(t) h_{\pm})(y)$. Then $u(y,t) = ([S_{+}(t) -
  S_{-}(t)] \half \ii B^{-1} h)(y)$. The operators $S_{\pm}(t)$ are
Fourier integral operators. Their construction is well known, see for
example Duistermaat \cite[Chapter 5]{duistermaat:1996}. Microlocally,
the solution operator associated with (\ref{eq:first_order_time}) can
be written in the matrix form
\[
   S(t) = V \left(\begin{array}{cc}
            S_{+}(t) & 0 \\ 0 & S_{-}(t)
            \end{array}\right) \Lambda ;
\]
in this notation, $S_{12}(t) = [S_{+}(t) - S_{-}(t)] \half \ii
B^{-1}$.

For the later analysis, we introduce the operators $S(t,s)$ and
$S_{\pm}(t,s)$: $S(t,s)$ solves the problem
\begin{eqnarray*}
   \left[\pdpd{}{t}
     - \left( \ba{cc} 0 & 1 \\ - A(x,D_x) & 0 \ea \right)
   \right] S(t,s) &=& 0 ,
\\
   S(\cdot,s) |_{t=s} &=& 0 ,\quad
   \pdpd{S}{t}(\cdot,s) |_{t=s} = \mathrm{Id} ,
\end{eqnarray*}
so that the solution of
\begin{equation*}
   \left[\pdpd{{}^2}{t^2} + A(x,D_x)\right] u = f ,
\quad u(t<0) = 0 ,
\end{equation*}
is given by
\begin{equation*}
   u(y,t) = \int_0^t \mathsf{P}_1 S(t,s)
   \left( \!\! \ba{c} 0 \\ f(\cdot,s) \ea \!\! \right)(y) \, \dd s
       = \iint G(y,x,t-s) f(x,s) \, \dd x \dd s ,
\end{equation*}
where we identified the causal Green's function $G(y,x,t-s)$. Here,
$\mathsf{P}_1$ is the projection, $\mathsf{P}_1
\left(\mstrut{0.32cm}\right. \!\!\!{\scriptsize \ba{cc} u
  \\ \pdpd{u}{t} \ea}\!\!\! \left.\mstrut{0.32cm}\right) =
u$. Likewise, $S_{+}(t,s)$ solves (for $t \in \R$) the problem
\begin{eqnarray*}
   P_{+}(x,D_x,D_t) \, S_{+}(\cdot,s) &=& 0 ,
\\
   S_{+}(\cdot,s) |_{t=s} &=& \mathrm{Id} ,
\end{eqnarray*}
so that the causal solution of
\begin{equation*}
   P_{+}(x,D_x,D_t) \, u_{+} = f_{+} ,\quad
   f_{+} = \half \ii B(x,D_x)^{-1} f ,
\end{equation*}
is given by
\begin{equation*}
   u_{+}(y,t) = \int_{-\infty}^t \!\! (S_{+}(t,s)
          f_{+}(\cdot,s))(y) \, \dd s
   = \iint G_{+}(y,x,t-s) f_{+}(x,s) \, \dd x \dd s ,
\end{equation*}
while the anticausal solution is given by
\begin{equation*}
   u_{+}(y,t) = - \! \int_t^{\infty} \! (S_{+}(t,s)
          f_{+}(\cdot,s))(y) \, \dd s
   = \iint G_{+}(y,x,s-t) f_{+}(x,s) \, \dd x \dd s .
\end{equation*}
A similar construction holds with $+$ replaced by $-$.

\subsection{Oscillatory integral representation}

For sufficiently small $t$ (in the absence of conjugate points), one
obtains the oscillatory integral representation,
\begin{equation}
\label{eq:smooth_Greens}
   (S_{\pm}(t) h_{\pm})(y)
   = (2\pi)^{-n}
     \iint a_{\pm}(y,t,\xi) 
     \exp[\ii \phi_{\pm}(y,t,x,\xi)] \,
     h_{\pm}(x) \, \dd x \dd \xi ,
\end{equation}
where
\begin{equation}
\label{eq:Maslov_phase}
   \phi_{\pm}(y,t,x,\xi)
   = \alpha_{\pm}(y,t,\xi) - \langle \xi,x \rangle .
\end{equation}
We note that $\alpha_{-}(y,t,\xi) = -
\alpha_{+}(y,t,-\xi)$. Singularities are propagated along the
bicharacteristics, which are determined by Hamilton's equations
generated by the principal symbol $\pm B^{\rm prin}(x,\xi)$
\begin{equation}
\label{eq:Hamilton}
   \pdpdo{y^t}{t} = \pm \pdpd{B^{\rm prin}(y^t,\eta^t)}{\eta} ,
\quad
   \pdpdo{\eta^t}{t} = \mp \pdpd{B^{\rm prin}(y^t,\eta^t)}{y} .
\end{equation}
We denote the solution of (\ref{eq:Hamilton}) with the $+$ sign and
initial values $(x,\xi)$ at $t = 0$ by $(y^t(x,\xi), \eta^t(x,\xi)) =
\Phi^t(x,\xi)$. The solution with the $-$ sign is found upon reversing
the time direction and is given by
$(y^{-t}(x,\xi),\eta^{-t}(x,\xi))$. Away from conjugate points, $y^t$
and $\xi$ determine $\eta^t$ and $x$; we write $x =
\widetilde{x}^t(y,\xi)$ and $\eta^t = \widetilde{\eta}^t(y,\xi)$. (We
also use the parametrization in which the roles of $(y,\xi)$ and
$(x,\eta)$ are interchanged.) Then
\[
   \alpha_{+}(y,t,\xi) = \langle \xi,\widetilde{x}^t(y,\xi) \rangle .
\]
To highest order,
\begin{equation}
\label{eq:amplitude_A_M2}
   a_{+}(y,t,\xi) = \left|\mstrut{0.5cm}\right.
       \left.\pdpd{(y^t)}{(x)}\right|_{x=\widetilde{x}^t(y,\xi),\xi} 
       \left.\mstrut{0.5cm}\right|^{-1/2} .
\end{equation}
We consider the perturbations of $(y^t,\eta^t)$ with respect to the
initial conditions $(x,\xi)$,
\begin{equation} \label{equ:fundmat}
W^t(x,\xi) = \left(
\begin{array}{ c c }
     W_{1}^t(x,\xi) & W_{2}^t(x,\xi) \\
     W_{3}^t(x,\xi) & W_{4}^t(x,\xi)
\end{array} \right)
= \left(
\begin{array}{ c c }
     \partial_x y^t(x,\xi) & \partial_\xi y^t(x,\xi) \\
     \partial_x \eta^t(x,\xi) & \partial_\xi \eta^t(x,\xi) 
\end{array} \right).
\end{equation} 
This matrix solves the (linearized) Hamilton-Jacobi equations,
\begin{equation} \label{equ:perturbS}
\pdpdo{W^t}{t}(x,\xi) = \left(
\begin{array}{rr}
     \partial_{\eta y} B^{\rm prin}(y^t,\eta^t) &
     \partial_{\eta\eta} B^{\rm prin}(y^t,\eta^t) \\
     -\partial_{yy} B^{\rm prin}(y^t,\eta^t) &
     -\partial_{y\eta} B^{\rm prin}(y^t,\eta^t) 
\end{array} \right) W^t(x,\xi) ,
\end{equation}
subject to initial conditions $W^{t=0} = \operatorname{I}$. We note
that away from conjugate points, the submatrix $W_{1}^t$ is
invertible. Because
\[
   \widetilde{x}^t = \pdpd{\alpha_{+}}{\xi} ,\quad
   \widetilde{\eta}^t = \pdpd{\alpha_{+}}{y} ,
\]
integration of \eqref{equ:perturbS} along $(y^t,\eta^t)$ yields:
\begin{eqnarray}
   \pdpd{^2 \alpha_{+}}{y\partial\xi}(y^t(x,\xi),t,\xi) &=&
   (W_{1}^t(x,\xi))^{-1} ,
\\
   \pdpd{^2 \alpha_{+}}{\xi^2}(y^t(x,\xi),t,\xi) &=&
   (W_{1}^t(x,\xi))^{-1} W_{2}^t(x,\xi) ,
\\
   \pdpd{^2 \alpha_{+}}{y^2}(y^t(x,\xi),t,\xi) &=&
   W_{3}^t(x,\xi) (W_{1}^t(x,\xi))^{-1} ,
\end{eqnarray}
which we evaluate at $x = \widetilde{x}^t(y,\xi)$. It follows that
\begin{equation*}
   a_{+}(y,t,\xi) = \left|\mstrut{0.3cm}\right.
       \det W^t_{1} |_{x=\widetilde{x}^t(y,\xi),\xi} 
       \left.\mstrut{0.3cm}\right|^{-1/2} .
\end{equation*}
The amplitude of $S_{+}(t) \, \half \ii B^{-1}$, then becomes
\[
   a_{+}(y,t,\xi) \, \half \ii
        B^{\mathrm{prin}}(\widetilde{x}^t(y,\xi),\xi)^{-1}
\]
to leading order; we denote this amplitude by
$\widetilde{a}_{+}(y,t,\xi)$. The amplitude $a_{-}$ follows from time
reversal: $a_{-}(y,t,\xi) = \overline{a_{+}(y,t,-\xi)}$.

In the case of conjugate points, we use the semigroup property of
$S(t,s)$ and decompose the time step into smaller time steps such that
in each step the formation of caustics is avoided.
Numerically, the size of the smaller time steps can be determined by monitoring the rank-deficiency of $W_1^t$, see \cite{dHUVW:2013} for a more general point of view and Subsection \ref{sec:comprtc} for an application.

\subsection{The source field}
\label{sec:source}
In the absence of caustics, we can change phase variables in the
oscillatory integral representation of $G$ according to
\begin{multline}
   G_{+}(y,x,t)
   = (2\pi)^{-1} \int
   \int
   (2\pi)^{-n} \int a_{+}(y,t',\xi) 
\\
     \exp[\ii \phi_{+}(y,t',x,\xi)]
     \, \dd \xi \, \exp[\ii \tau (t - t')] \,
      \dd t' \dd\tau
\\
   = (2\pi)^{-1}
     \int a_{+}'(y,x,\tau)
     \exp[\ii \tau (t - T(y,x))] \, \dd\tau .
\end{multline}
By applying the method of stationary phase in the variables $(\xi,t')$, one can show that the source field can be written in the form  \cite{Brytik2013}
\begin{equation}
\label{eq:GMOI}
   G(x,\tilde{x},t)
   = (2\pi)^{-1}
     \int a'(x,\tilde{x},\tau)
     \exp[\ii \tau (t - T(x,\tilde{x}))] \, \dd\tau .
\end{equation}
Here, $\tilde{x}$ is the source location and $T$ is the travel time
satisfying the eikonal equation
\begin{equation}
\label{eq:eik}
   B^{\mathrm{prin}}(x,-\partial_x T(x,\tilde{x})) = -1
\end{equation}
and $a' = \mathcal{A}$ to highest order with
\begin{equation}
   | \mathcal{A}(x,\tilde{x},\tau) | =
(2\pi)^{-(n-1)/2}
   \left|\mstrut{0.45cm}\right.
   \det \pdpd{(x,\xi,t)}{(y,x,\tau)}
        \left.\mstrut{0.45cm}\right|^{1/2}
\end{equation}
see \cite{Brytik2013} for details. We introduce
\begin{equation}
\label{eq:sl}
   n_{\tilde{x}}(x) = \frac{\partial_x T(x,\tilde{x})}{
                      |\partial_x T(x,\tilde{x})|} ;
\end{equation}
in view of (\ref{eq:eik}),
\[
   |\partial_x T(x,\tilde{x})| =
      \frac{1}{B^{\mathrm{prin}}(x,n_{\tilde{x}}(x))} .
\]
We note that through $n_{\tilde{x}}(x)$ we obtain the incidence angle
of the source field at $x$.  In Section \ref{sec:angle}, we will
arrange and study the images with respect to incidence angle.  We also
note that $n_{\tilde{x}}(x)$ can be estimated from the Poynting vector
$ -\partial_t G(x,\tilde{x},t) \ \partial_x G(x,\tilde{x},t) $ at $t =
T(x,\tilde{x})$ \cite{Yoon-Marfurt-Starr_2004,Yoon-Marfurt_2006} or
from $ -\partial_t G(x,\tilde{x},t) \stackrel{(t=0)} {\ast}\ \partial_x G(x,\tilde{x},-t) , $ 
(possibly normalized by the autocorrelation, $\partial_t G(x,\tilde{x},t)  \stackrel{(t=0)} {\ast} \partial_t G(x,\tilde{x},-t)$; note that this normalization is primarily applied to suppress the dependency on $a'$), for instance in the PDE solution formulation of RTM.

\section{Reverse-time continuation from the boundary}
\label{sec:3}
The key results we obtain in this section are the formulation of an oscillatory integral representation and its computation using dyadic parabolic decomposition and wave packets for reverse-time continuation with a boundary source. These are also central in the formulation and computation of the inverse scattering and imaging operators presented in Section \ref{sec:4}.
We introduce Euclidean boundary normal coordinates, $x = (x',x_n)$;
that is, $x' = (x_1,\ldots,x_{n-1})$, and $x_n = 0$ defines the
boundary. We let $\Sigma$ denote a bounded open subset of $\{ (x,t)
\in \R^n_x \times \R_t\ |\ x_n = 0 \}$. We denote the restriction to
the boundary by $R_{x_n}$.

We let $w_r$ be an anticausal solution to
\begin{equation} \label{eq:revtcont}
   \left[\pdpd{{}^2}{t^2} + A(x,D_x)\right] \, w_r(x,t)
    = \delta(x_n) \, g(x',t) ;
\end{equation}
we have $w_r = w_{r,+} + w_{r,-}$ with
\begin{equation*}
   w_{r,+}(y,t) = - \! \int_t^{\infty} \! (S_{+}(t,s)
      \half \ii B^{-1} R_{x_n}^*
      \widetilde{\mathit{\Psi}}_{\Sigma}
                  g(\cdot,s))(y) \, \dd s 
\end{equation*}
noting that
\[
   R_{x_n}^* g(x,t) = \delta(x_n) \, g(x',t)
\]
for any functions $g$ defined on $\R^{n-1}_{x'} \times \R_t$. Here,
$\widetilde{\mathit{\Psi}}_{\Sigma} =
\widetilde{\mathit{\Psi}}_{\Sigma}(x',t,D_{x'},D_t)$ is a
pseudo{\-}differential cutoff designed to remove grazing rays. The
relation between contributions from negative frequencies and positive
frequencies is
\begin{equation} \label{eq:pm}
   w_{r,-}(y,t) = \overline{w_{r,+}(y,t)} .
\end{equation}

We now introduce principal parts of symbols,
$C_{\pm}(x',x_n,\xi',\tau)$, as the solutions for $\zeta$ of
\[
   A^{\mathrm{prin}}(x',x_n,\xi',\zeta) = \tau^2 .
\]
We write $C(x',x_n,\xi',\tau) = C_+(x',x_n,\xi',\tau)$. In the further
analysis we will need the operator,
\[
   C(x',x_n,D_t^{-1} D_{x'},1)\quad\text{at the surface, $x_n = 0$}
\]
with principal symbol $C(x',x_n,\tau^{-1} \xi',1)$.


\subsection{Conjugate points}

In the case of conjugate points, we introduce a partition of unity into
$\Sigma \subset \R^{n-1}_{x'} \times \R_t$ (with overlap $\delta$ in
time). Incorporating this partition of unity in
$\widetilde{\mathit{\Psi}}_{\Sigma}$, we obtain a set of cutoffs,
$\widetilde{\mathit{\Psi}}_{\Sigma,ij}$. The first index signifies a
subdivision in $\R^{n-1}_{x'}$ while the second index identifies
intervals in time.

To describe the use of the semigroup property, we fix $i$. Let $i=1$
and assume, without loss of generality, that we need two smaller time intervals, $[t,t+t_1]$ and $[t + t_1,
  T_1]$, say, to avoid conjugate points in the smaller time intervals. Then
\begin{multline}
   \int_t^{T_1} \! (S_{+}(t,s)
     \half \ii B^{-1} R_{x_n}^*
     (\widetilde{\mathit{\Psi}}_{\Sigma,11}
        + \widetilde{\mathit{\Psi}}_{\Sigma,12})
                  g(\cdot,s))(y) \, \dd s
\\
   = \int_t^{t+t_1} \! (S_{+}(t,s)
     \half \ii B^{-1} R_{x_n}^*
     \widetilde{\mathit{\Psi}}_{\Sigma,11}
                  g(\cdot,s))(y) \, \dd s
\\
   + S_{+}(t,t+t_1-\delta) \int_{t+t_1-\delta}^{T_1} \!
      (S_{+}(t+t_1-\delta,s) \half \ii B^{-1} R_{x_n}^*
     \widetilde{\mathit{\Psi}}_{\Sigma,12}
                  g(\cdot,s))(y) \, \dd s .
\end{multline}
We now focus on representations for
\begin{multline*}
   \int_t^{t+t_1} \! (S_{+}(t,s)
     \half \ii B^{-1} R_{x_n}^*
     \widetilde{\mathit{\Psi}}_{\Sigma,11}
                  g(\cdot,s))(y) \, \dd s
\quad\text{and}
\\
   \int_{t+t_1-\delta}^{T_1} \!
      (S_{+}(t+t_1-\delta,s) \half \ii B^{-1} R_{x_n}^*
     \widetilde{\mathit{\Psi}}_{\Sigma,12}
                  g(\cdot,s))(y) \, \dd s ,
\end{multline*}
in the absence of conjugate points.

\subsection{Oscillatory integral representations}

We have
\begin{multline}
   \chi_n \, w_{r,+}(y,t) = (2\pi)^{-n}
   \int \!\! \int_t^{\infty} \!\! \int
      \chi_n \, a^{(\mathrm{bkd})}(x',s - t,y,\eta) \,
\\
   \exp[\ii (-\alpha_{+}(x',0,s - t,\eta)
        + \langle \eta,y \rangle)] \,
   g(x',s) \, \dd x' \dd s \, \dd\eta ,
\end{multline}
where
\begin{equation}
   a^{(\mathrm{bkd})}(x',s - t,y,\eta) =
   \left|\mstrut{0.5cm}\right. \!
   \left.\pdpd{(y^{s-t})}{(x)}\right|_{\eta,x=x^{s-t}(x',0,\eta)} 
       \! \left.\mstrut{0.5cm}\right|^{-1/2}
   \tfrac{1}{2} \ii \tau^{-1}
       \widetilde{\mathit{\Psi}}_{\Sigma}(x',s,\xi',\tau)
\end{equation}
up to terms of lower order, that is, the error (expressed in $\eta$) is of order $(1 + |\eta|^2)^{-1}$, and
$\chi_n$ is a cutoff function which removes contributions for $x_n < 0$ 
(the expressions for $\xi'$ and $\tau$ in terms of $\eta$ are given in \eqref{eq:dalp2} and \eqref{eq:dalp3} below).
The operator $\chi_n \, S_{+}(t,s) \half \ii
B^{-1} R_{x_n}^* \widetilde{\mathit{\Psi}}_{\Sigma}$ is a FIO, the canonical relation of which is a subset of
\[
   \{ (y,\eta;(y^{s-t})'(y,\eta),s-t,(\eta^{s-t})'(y,\eta),
                    -B^{\mathrm{prin}}(y,\eta))\ |\
            (y^{s-t})_n(y,\eta) = 0 \} .
\]
The dyadic parabolic decomposition of phase space enters in the
reverse-time continuation as
\begin{multline}
   \chi_n(y_n) \, w_{r,+}(y,t) = \chi_n(y_n) \, \sum_{\nu,k}
   \iint
   \left\{ \mstrut{0.5cm} \right. (2\pi)^{-n}
   \int \widehat{\beta}_{\nu,k}(\eta) \widehat{\chi}_{\nu,k}(\eta)
\\
      a^{(\mathrm{bkd})}(x',s - t,y,\eta) \,
   \exp[-\ii \alpha_{+}(x',0,s - t,\eta)] \, \dd\eta
   \left. \mstrut{0.5cm} \right\}
          g(x',s) \, \dd x' \dd s .
\end{multline}
Fixing $(\nu,k)$ corresponds with (directional) ``controlled
illumination.''

\subsection{Boundary source decomposition; wave packets in space-time}
\label{sec:33}
We change phase variables in the representation for $w_{r,+}$. We
could do this in two steps, changing parametrizations from
$((x',x_n),\eta)$ to $(y,(\xi',\xi_n))$ and then $(s,\xi_n)$ to
$(x_n,\tau)$. Here, we carry out this change in a single step:
\begin{multline}
   \chi_n \, w_{r,+}(y,t) = (2\pi)^{-2n}
   \iint \!\! \int_t^{\infty} \!\! \iint
      a^{(\mathrm{bkd})}(x',s - t,y,\eta) \,
\\
   \exp[\ii (-\alpha_{+}(x',0,s - t,\eta)
        + \langle \eta,y \rangle)] \,
   \exp[\ii (\tau \, s
            + \langle \xi',x' \rangle)] \,
   \dd\eta \, \dd x' \dd s\
   \widehat{g}(\xi',\tau) \, \dd\xi' \dd\tau ;
\end{multline}
applying the method of stationary phase in $(\eta,x',s)$ yields
solving
\begin{eqnarray}
   \partial_{\eta} \alpha_{+}(x',0,s - t,\eta) &=& y    , \label{eq:dalp1}
\\
     \partial_{x'} \alpha_{+}(x',0,s - t,\eta) &=& \xi' , \label{eq:dalp2}
\\
        \partial_s \alpha_{+}(x',0,s - t,\eta) &=& \tau \label{eq:dalp3}
\end{eqnarray}
for given $(y,\xi',\tau)$ and $t$ fixed (which is viewed as a
parameter here). The solutions, $(\eta_0,x'_0,s_0)$, are the
stationary points of $-\alpha_{+}(x',0,s - t,\eta) + \langle \eta,y
\rangle + \tau \, s + \langle \xi',x' \rangle$. We have $s_0 >
  t$. These equations imply that
\begin{equation*}
\left.\begin{array}{rcr}
   y \!\! &=& \!\! \widetilde{y}^{s_0-t}(x'_0,0,\eta_0)
\\[0.2cm]
\xi' \!\! &=& \!\! \widetilde{\xi}^{s_0-t}{\,}'(x_0',0,\eta_0)
\end{array}\right\}\ \text{that is,}\ \
   (x_0',0,\widetilde{\xi}^{s_0-t}{\,}',
   C(x_0',0,\widetilde{\xi}^{s_0-t}{\,}',\tau))
      \stackrel{\Phi^{s_0-t}}{\to}
          (\widetilde{y}^{s_0-t},\eta_0).
\end{equation*}
For given $t$, $s_0$ is determined since $(x_0',0)$, $\eta_0$ and $s_0 - t$
determine a unique ray, in view of the absence of conjugate
points. Thus we need to solve
\begin{eqnarray}
   \eta_0 &=&
    \,\, \widetilde{\eta}^{s_0-t}(y,\xi',C(x_0',0,\xi',\tau)) ,
\label{eq:statpts-1}
\\
     x'_0 &=& \widetilde{x}^{s_0-t}{\,}'(y,\xi',C(x_0',0,\xi',\tau)) ,
\label{eq:statpts-2}
\\
   0 &=& \ \widetilde{x}^{s_0-t}_n(y,\xi',C(x_0',0,\xi',\tau))\ \ \
   \text{or}\ \ s_0 = T(x_0',0,y) + t ,
\label{eq:statpts-3}
\end{eqnarray}
for $(\eta_0,x_0',s_0)$. To obtain a unique solution, in general, we
need to localize $g$, which we do by substituting a wave packet
contribution, that is, $g_{\gamma}
\widehat{\varphi}_{\gamma}(\xi',\tau)$ for
$\widehat{g}(\xi',\tau)$. Then
\begin{multline}
   -\alpha_{+}(x'_0(y,\xi',\tau;t),0,s_0(y,\xi',\tau;t) - t,
               \eta_0(y,\xi',\tau;t))
     + \langle \eta_0(y,\xi',\tau;t),y \rangle
\\
   = -\langle \eta_0(y,\xi',\tau;t),
     y^{s_0(y,\xi',\tau;t)-t}(x'_0(y,
        \xi',\tau;t),0,\eta_0(y,\xi',\tau;t)) \rangle
\\
              + \langle \eta_0(y,\xi',\tau;t),y \rangle = 0 
\end{multline}
while
\begin{multline}
   \tau \, s_0(y,\xi',\tau;t) + \langle \xi',x'_0(y,
        \xi',\tau;t) \rangle
\\
   = \tau \, s_0(y,\xi',\tau;t) + \langle \xi',
        \widetilde{x}^{s_0(y,\xi',\tau;t)-t}{\,}'(y,
           \xi',C(x_0'(y,\xi',\tau;t),0,\xi',\tau)) \rangle
\\
   =: \theta_+(y,t,\xi',\tau) .
\end{multline}
We evaluate
\begin{multline}
   \left.\pdpd{{}^2 [-\alpha_{+}(x',0,s - t,\eta)
        + \langle \eta,y \rangle + \tau \, s
            + \langle \xi',x' \rangle]}{\eta \partial\eta}
   \right|_{(\eta_0,x'_0,s_0)}
\\
   = -\left.\pdpd{{}^2 [\alpha_{+}(x',0,s - t,\eta)]}{\eta \partial\eta}
   \right|_{(\eta_0,x'_0,s_0)}
   = -(W^{s_0-t}_1(y,\eta_0))^{-1} W^{s_0-t}_2(y,\eta_0) ,
\end{multline}
\begin{multline}
   \left.\pdpd{{}^2 [-\alpha_{+}(x',0,s - t,\eta)
        + \langle \eta,y \rangle + \tau \, s
            + \langle \xi',x' \rangle]}{\eta \partial x'}
   \right|_{(\eta_0,x'_0,s_0)}
\\
   = -\left.\pdpd{{}^2 \alpha_{+}(x',0,s - t,\eta)}{\eta \partial x'}
   \right|_{(\eta_0,x'_0,s_0)} = -[(W^{s_0-t}_1(y,\eta_0))^{-1}]' ,
\end{multline}
and
\begin{multline}
   \left.\pdpd{{}^2 [-\alpha_{+}(x',0,s - t,\eta)
        + \langle \eta,y \rangle + \tau \, s
            + \langle \xi',x' \rangle]}{x' \partial x'}
   \right|_{(\eta_0,x'_0,s_0)}
\\
   = -\left.\pdpd{{}^2 \alpha_{+}(x',0,s - t,\eta)}{x' \partial x'}
   \right|_{(\eta_0,x'_0,s_0)} = - \, {}'[W^{s_0-t}_3(y,\eta_0)]\
            [(W^{s_0-t}_1(y,\eta_0))^{-1}]' ,
\end{multline}
subject to the substitutions according to
(\ref{eq:statpts-1})-(\ref{eq:statpts-3}), and then
\begin{multline}
   \left.\pdpd{{}^2 [-\alpha_{+}(x',0,s - t,\eta)
        + \langle \eta,y \rangle + \tau \, s
            + \langle \xi',x' \rangle]}{s^2}
   \right|_{(\eta_0,x'_0,s_0)}
\\
   = -\left.\pdpd{{}^2 \alpha_{+}(x',0,s - t,\eta)}{s^2}
   \right|_{(\eta_0,x'_0,s_0)}
                 = -\left.\pdpd{\tau}{s} \right|_{s = s_0 - t} ,
\end{multline}
\begin{multline}
   \left.\pdpd{{}^2 [-\alpha_{+}(x',0,s - t,\eta)
        + \langle \eta,y \rangle + \tau \, s
            + \langle \xi',x' \rangle]}{s \partial \eta}
   \right|_{(\eta_0,x'_0,s_0)}
\\
   = -\left.\pdpd{{}^2 \alpha_{+}(x',0,s - t,\eta)}{s \partial \eta}
   \right|_{(\eta_0,x'_0,s_0)}
         = -\left.\pdpd{\widetilde{y}^{s-t}}{s}
                         \right|_{s = s_0} \!\! (x_0',0,\eta_0) ,
\end{multline}
and
\begin{multline}
   \left.\pdpd{{}^2 [-\alpha_{+}(x',0,s - t,\eta)
        + \langle \eta,y \rangle + \tau \, s
            + \langle \xi',x' \rangle]}{s \partial x'}
   \right|_{(\eta_0,x'_0,s_0)}
\\
   = -\left.\pdpd{{}^2 \alpha_{+}(x',0,s - t,\eta)}{s \partial x'}
   \right|_{(\eta_0,x'_0,s_0)}
                 = -\left.\pdpd{{\widetilde{\xi}^{s-t}}{\,}'}{s}
                         \right|_{s = s_0} \!\! (x_0',0,\eta_0) .
\end{multline}
From these expressions we form the Hessian which is used to transform
$a^{(\mathrm{bkd})}(x',s - t,y,\eta)$ to the amplitude
$a_b^{(\mathrm{bkd})}(y,t,\xi',\tau)$, so that
\begin{multline} \label{eq:wprepr}
   \chi_n \, w_{r,+}(y,t) = (2\pi)^{-n} \sum_{\gamma}
           g_{\Sigma,\gamma} \iint \chi_n \,
   a_b^{(\mathrm{bkd})}(y,t,\xi',\tau)
\\
   \exp[\ii \, \theta_+(y,t,\xi',\tau)] \,
   \widehat{\varphi}_{\gamma}(\xi',\tau) \, \dd\xi' \dd\tau .
\end{multline}
Essentially, this representation corresponds with local coordinates
$(x_n,y,\xi',\tau)$ for the canonical relation of the solution
operator with $t$ fixed.

\subsection{Algorithm}
\label{sec:comprtc}

We adapt the ''box algorithm'' for the multiscale discrete approximation of FIOs developed in \cite{AndersondHW:2012} to \eqref{eq:wprepr}, with accuracy $\mathcal{O}(2^{-k/2})$ at frequency scale $k$. In the general case, the medium can be laterally varying at the boundary. Then we need to employ compactly supported cutoff functions in $x'$, realized by the partition of unity $\widetilde{\mathit{\Psi}}_{\Sigma,ij}$. Within each cutoff $i$ the lateral variation of the (smooth) velocity model is assumed to be negligible at the boundary, and the algorithm outlined below can then be applied for each cutoff $i$ separately.

For convenience of notation, we now assume that the wave speed does not vary laterally at the surface and fix $i=1$. 
Without loss of generality, we suppose that we need $N_s$ time intervals $[t+(n_s-1)t_1,t+n_st_1]$, $n_s=1,\dots,N_s$, of duration $t_1=(T_1-t)/N_s$ in order to avoid the formation of conjugate points. 
Numerically, such a covering of sub-time intervals can be determined straightforwardly from the points of rank-deficiency of the matrix $W_1^t$ of the Hamiltonian system as detailed in \cite{dHUVW:2013} and sketched below.

Let us consider one frequency box $\hat\chi_{\nu,k}$. We begin with computing the bicharacteristics (rays) of the Hamiltonian system, $(x_0',0,\nu)\mapsto (y,\eta)=(\tilde y^s,\tilde \eta^s)$, i.e. $(x_0',0,\nu)=\pdpd{\theta_+(y,t,\xi',\tau)}{(\xi',\tau)}$ where $s\in(0,t_1]$. For each time interval $n_s$, we thus obtain the \emph{coordinate transform} \cite{AndersondHW:2012}
$$
T^{(n_s)}_{\nu,k}(y)=\left(x_0',s+t+(n_s-1)t_1\right).
$$
The solution of the corresponding Hamilton-Jacobi system yields the propagator matrix $W^s$ from which we obtain the quantities
\begin{eqnarray*}
\pdpd{^2\theta_+(y,t,\xi',\tau)}{y\partial(\xi',\tau)} &=& \left(W_1^s\right)^{-1}\\
\pdpd{^2\theta_+(y,t,\xi',\tau)}{^2(\xi',\tau)} &=& -\left(W_1^s\right)^{-1}W_2^s\\
\pdpd{^2\theta_+(y,t,\xi',\tau)}{^2y} &=& W_3^s\left(W_1^s\right)^{-1}.
\end{eqnarray*}
We can now apply the box algorithm to each time interval $n_s$ and obtain the (partially) reverse-time continued wave field from (data) boundary sources
\begin{multline}
\label{equ:rtcbox}
w_{r,+}^{(n_s)}(y,t+(n_s-1)t_1)= 
\sum_{\nu,k}a^{(\textnormal{bkd})}(y,\nu)\sum_{r=1}^{R_{\nu,k}} \alpha_{\nu,k}^{(r)}(y)\\
\sum_{(\xi',\tau)} e^{\ii \langle \Tvk^{(n_s)}(y),(\xi',\tau)\rangle}\hat g(\xi',\tau)\hat\beta_{\nu,k}(\xi',\tau) \hat\chi_{\nu,k}(\xi',\tau)\hat\vartheta_{\nu,k}^{(r)}(\xi',\tau)
\end{multline}
where $\alpha_{\nu,k}^{(r)}$ and $\vartheta_{\nu,k}^{(r)}$ are the expansion functions arising in the tensor-product representation of the complex exponential of the second-order Taylor expansion term of $\theta_+$ on the frequency box $\hat\chi_{\nu,k}$ \cite{AndersondHW:2012}.

\begin{algorithm}
\tt
\begin{description}
\item[\sc Part I -- reverse-time continuation from the boundary, semigroup 1] \hfill
\begin{description}
\item[{\rm for $n_s=1:N_s$}] \hfill\\ 
1. compute coordinate transforms $\Tvk$ and propagator matrices $W$\\
2. compute $w_{r,+}^{(n_s)}(y,t+(n_s-1)t_1)$: box algorithm,  \eqref{equ:rtcbox}
\item[{\rm end}]
\end{description}
\item[\sc Part II -- half wave equation reverse-time continuation, semigroups]\hfill
\begin{description}
\item[{\rm for $n_p=2:N_s$}]\hfill
\begin{description}
\item[{\rm for $n_s=n_p:N_s$}] \hfill\\ 
half wave equation evolution operator $P_{t_1}$: box algorithm\\
\qquad$w_{r,+}^{(n_s)}(y,t+(n_s-n_p)t_1)=P_{t_1}w_{r,+}^{(n_s)}(y,t+(n_s-n_p+1)t_1)$
\item[{\rm end}]
\end{description}
\item[{\rm end}]
\end{description}
\end{description}
{\sc Wave Field}\quad$w_{r,+}(y,t)=\sum_{n_s=1}^{N_s}w_{r,+}^{(n_s)}(y,t)$
\caption{\label{tab:rtc} Outline of reverse-time continuation from the boundary in the case of conjugate points. In the absence of caustics, the algorithm reduces to Part I, with $N_s=1$.}
\end{algorithm}

To obtain the final reverse-time continued wave fields $w_{r,+}^{(n_s)}(y,t)$, we construct a parametrix for the Cauchy initial value problems for the half wave equation with initial data $w_{r,+}^{(n_s)}(y,t+(n_s-1)t_1)$, $n_s=2,\dots,N_s$, initial time $t+(n_s-1)t_1$ and final time $t$. We compute these parametrices using the box algorithm (this has been studied in detail in \cite{AndersondHW:2012}). We make use of the semigroup property and obtain the parametrix for the reverse-time interval $[t+(n_s-1)t_1, t]$ as the composition of the parametrices for the time intervals $[t+(n_s-n_p+1)t_1,t+(n_s-n_p)t_1]$, $n_p=2,\dots,n_s$.
Finally, we have
$$
w_{r,+}(y,t)=\sum_{n_s=1}^{N_s}w_{r,+}^{(n_s)}(y,t).
$$
The different steps involved in modeling receiver wave propagation from the boundary in reverse-time are summarized in Algo. \ref{tab:rtc} and illustrated in Fig. \ref{FigBPalgo} for a numerical example that is detailed in Section \ref{sec:5}.

The coordinate transform $T^{(n_s)}_{\nu,k}(y)$ and the propagator matrix $W^s$ can numerically be evaluated as follows. 
For simplicity, we consider the case of isotropic medium.
Let $c$ be the wave speed at the boundary and $\nu=(\nu',\nu_n)=(\xi',\tau)/|(\xi',\tau)|$. Then $T^{(n_s)}_{\nu,k}(y)$ follows from the bicharacteristics (rays) of the half wave equation supplemented with initial conditions $y^0=(x',0)$, $\eta^0=({\eta^0}',\eta_n^0)=\frac{c}{\nu_n}(\nu',\sqrt{\nu_n^2/c^2-|\nu'|^2})$ (for evolution time $s$). Similarly, $W^s$ is obtained as the solution of the Hamilton-Jacobi system associated with the half wave equation with initial conditions $W_2^0=W_3^0\equiv0$, $W_1^0=\left(\begin{array}{cc}\mathbf{I}_{n-1}&0\\c\eta'&c\eta_n\end{array}\right)$ and $W_4^0=\frac{c}{\nu_n}\left(\begin{array}{cc}\mathbf{I}_{n-1}&\frac{c\nu'}{\nu_n\eta_n}\\\frac{\nu'}{\nu_n}&\frac{c |\nu'|^2}{\nu_n^2\eta_n}\end{array}\right)$.

Finally, the duration $t_1$ for the time intervals is fixed numerically to be smaller than the length of the largest time interval $(0, t^*]$ for which $W_1^s$, $s\in(0,t^*]$, is nonsingular for the discrete set of values for $\nu$ considered. Note that if conjugate points are detected in the subsequent time-stepping intervals, the concerned time intervals can be further broken up into intervals of smaller size, again using the semigroup property, without the need to recompute the reverse-time continuation up to these points.

\section{Inverse scattering}
\label{sec:4}

We assume that a source at $\tilde{x}$ generates the data,
$d_{\Sigma}(x',t)$. We introduce the pseudodifferential operator \cite{Brytik2013}
\begin{equation}
\label{eq:Ndef}
   \mathcal{N}(x',D_{x'},D_t) =
   -2 \ii D_t \pdpd{B^{\rm prin}}{\xi_n}(x',0,D_t^{-1} D_{x'},
          C(x',0,D_t^{-1} D_{x'},1)).
\end{equation}
Furthermore, we introduce the pseudodifferential cutoff,
$\mathit{\Psi}_{\Sigma}$, which acts as a smooth cutoff which goes to
zero near $\partial\Sigma$, removes direct rays, and removes grazing
rays; that is, its symbol vanishes where
\[
   \pdpd{B^{\rm prin}}{\xi_n}(x',0,\tau^{-1} \xi',
          C^{\rm prin}(0,x',\tau^{-1} \xi',1)) = 0.
\]
These cutoffs commute up to leading order (through the product of their symbols), $\widetilde{\mathit{\Psi}}_{\Sigma} \mathit{\Psi}_{\Sigma}  = {\mathit{\Psi}}_{\Sigma} \widetilde{\mathit{\Psi}}_{\Sigma}$, which follows from the standard calculus of pseudodifferential operators \cite{Taylor1981}.

We let $w_r$ be an anticausal solution of (\ref{eq:revtcont}) with
\begin{equation}
   g(x',t) = \mathcal{N}(x',D_{x'},D_t) \,
      \mathit{\Psi}_{\Sigma}(x',t,D_{x'},D_t)
                   d_{\Sigma}(x',t).
\end{equation}
We define first-order partial differential and pseudodifferential
operators $\mathit{\Xi}(x,D_x,D_t)$ and $\mathit{\Theta}(x,D_x,D_t)$
with (principal) symbols
\begin{eqnarray*}
   \mathit{\Xi}_0(x,\xi,\tau) = \tau,&\quad &
   \mathit{\Xi}_j(x,\xi,\tau) = \xi_j \\
   \mathit{\Theta}_0(x,\xi,\tau) =
   \tau ,&\quad &
   \mathit{\Theta}_j(x,\xi,\tau) = \tau \,
         \pdpd{B^{\rm prin}}{\xi_j}(x,\xi) .
\end{eqnarray*}
We then define the pseudodifferential operator $L$ and the operator $K$ as
\begin{equation}
\label{equ_ISV_G}
\begin{split}
   L w(y,t) &= \mathcal{A}(y,\tilde{x},D_t)^{-1}
         2\ii D_t \sum_{p=0}^n
   \mathit{\Xi}_p(y,-\partial_y T(y,\tilde{x}),1)
   \mathit{\Theta}_p(y,D_y,D_t)
   w(y,t) ,
\\
 K w(y)   &= w(y,T(y,\tilde{x})).
   \end{split}
\end{equation}
%
%
Operator $K$ is a restriction to a hypersurface in
$\R^{n+1}$. 
The imaging operator, $H$, is then defined as
\[
   H d_{\Sigma}(y) = (K L (w_{r,+} + w_{r,-}))(y) .
\]
To leading order symbols, we get
\begin{multline}
\label{equ:is:final}
   L \chi_n w_{r,+}(y,t) = \frac{1}{(2\pi)^n}
\mathcal{A}(y,\tilde{x},D_t)^{-1} 2\ii D_t
   \iint  \chi_n
   a_b^{(\mathrm{bkd})}(y,t,\xi',\tau)
\\
\times \sum_{p=0}^n
   \mathit{\Xi}_p(y,-\partial_y T(y,\tilde{x}),1)
   \mathit{\Theta}_p(y,\partial_y \theta_{+},
                       \partial_t \theta_{+}) 
   w(y,T(y,\tilde{x}))
\\[0.2cm]
\times   \exp[\ii \, \theta_+(y,t,\xi',\tau)] \,
   \widehat{g}(\xi',\tau) \, \dd\xi' \dd\tau .
\end{multline}
%

\subsection{Isotropic case}

In the isotropic case, 
\begin{eqnarray*}
   A^{\mathrm{prin}}(x,\xi) &=& c(x)^2 \xi^2 ,
\\[0.2cm]
   B^{\mathrm{prin}}(x,\xi) &=& c(x) \, |\xi| ,
\\[0.2cm]
   C(x',0,\tau^{-1} \xi',1) &=&
               \sqrt{1 - c(x',0)^2 \tau^{-2} {\xi'}^2} ,
\\
   \mathit{\Theta}_0(x,\xi,\tau) &=& \tau ,\qquad
   \mathit{\Theta}_j(x,\xi,\tau) = \tau c(x) \frac{\xi_j}{|\xi|}
\end{eqnarray*}
and \eqref{equ:is:final} simplifies to
\begin{multline}
   L \chi_n w_{r,+}(y,t) = \frac{1}{(2\pi)^n}\frac{1}{\mathcal{A}_{\textrm{g}}(y,\tilde{x})}
    \partial_t^{-\!\frac{n+1}{2}}
   \int 
    (\ii\tau)^{(n-3)/2}
    \int \chi_n
   a_b^{(\mathrm{bkd})}(y,t,\xi',\tau)
\\
\times   \ii \, [\partial_t \theta_+(y,t,\xi',\tau)
      + c(y) n_{\tilde{x}}(y) \cdot
        \partial_y \theta_+(y,t,\xi',\tau)]
\\[0.2cm]
\times    \exp[\ii \, \theta_+(y,t,\xi',\tau)] \,
   \widehat{g}(\xi',\tau) \, \dd\xi' \dd\tau ,
\end{multline}
using that
\[
\mathcal{A}(y,\tilde{x},\tau)=\mathcal{A}_{\textrm{g}}(y,\tilde{x})(\ii\tau)^{(n-3)/2}.
\]
Operator $\partial_t^{-\!\frac{n+1}{2}}$ is to be read as the pseudodifferential operator with symbol $\tau \mapsto \tilde\sigma(\tau)(\ii \tau)^{-\!\frac{n+1}{2}}$ in which $\tilde\sigma$ is a smooth function, valued 1 except for the origin where it is 0.

\begin{algorithm}
\begin{description}
\tt
\item[\sc Part I -- boundary reverse-time continuation partial image] \hfill
\begin{description}
\item[{\rm for $n_s=1:N_s$}] \hfill\\ 
1. compute coordinate transforms $\Tvk$ and propagator matrices $W$\\
2. compute $w_{r,+}^{(n_s)}(y,t+(n_s-1)t_1)$: box algorithm,  \eqref{equ:rtcbox}\\
3. determine image region $y^*$, coordinate transform $\Tvk^{*(n_s)}(y^*)$, \\
\textcolor{white}{.}\qquad propagator matrices $W^{s^*}$\\
4. compute partial image $\Delta^{*(n_s,t+(n_s-1)t_1)}(y)$: box algorithm,  \eqref{equ:rtcimbox}
\item[{\rm end}]
\end{description}
\item[\sc Part II -- half wave equation reverse-time continuation partial image]\hfill
\begin{description}
\item[{\rm for $n_p=2:N_s$}]\hfill
\begin{description}
\item[{\rm for $n_s=n_p:N_s$}] \hfill\\ 
1. half wave equation evolution operator $P_{t_1}$: box algorithm\\
\qquad$w_{r,+}^{(n_s)}(y,t+(n_s-n_p)t_1)=P_{t_1}w_{r,+}^{(n_s)}(y,t+(n_s-n_p+1)t_1)$\\
2. image region $\tilde y^{s*}$: $s^*=t+(n_s-n_p+1)t_1-T(\tilde y^{s*},\tilde{x})$\\
3. coordinate transform, propagator matrices\hfill\\
\qquad $\tilde T_{\nu,k}^{*(n_s)}(\tilde y^{s*})=y_0$, $\tilde W^*=W(y_0,\eta_0,s^*)$\\
4. evaluate partial image $\tilde \Delta^{*(n_s,n_p)}(y)$: box algorithm
\item[{\rm end}]
\end{description}
\item[{\rm end}]
\end{description}
\end{description}
{\sc Image}\quad $\Delta_{d_\Sigma}(y)=\sum_{n_s=1}^{N_s} \Delta^{*(n_s,t+(n_s-1)t_1)}(y) + \sum_{n_p=2}^{N_s}\sum_{n_s=n_p}^{N_s}\tilde \Delta^{*(n_s,n_p)}(y)$
\caption{\label{tab:rtcim} Outline of inverse scattering in the case of conjugate points. In the absence of caustics, the algorithm reduces to Part I, with $N_s=1$.}
\end{algorithm}

\subsection{Computation}
We can use (\ref{eq:pm}) in the computations.
Through a simple modification, we can incorporate the imaging condition in the box algorithm for reverse-time continuation from the boundary detailed in Section \ref{sec:comprtc}, yielding an RTM imaging algorithm. Without loss of generality, we assume here that the source signature is a delta function; general discrete source signatures can be accommodated for in a straight-forward way by viewing them as a weighted sum of delta functions shifted by the time step for solving the Hamilton-Jacobi equations.
Suppose that the source travel time $T(y,\tilde{x})$ and amplitude $\mathcal{A}(y,\tilde{x},\tau)$ have been evaluated for the image region (here, by evaluating the corresponding Hamiltonian and Hamilton-Jacobi system, i.e. ray-tracing; cf. Section \ref{sec:source}). We begin with the evaluation of the imaging operator $H$ for partial reverse-time continuation from the boundary (cf. Algo. \ref{tab:rtc}, Part I). We obtain a contribution of time interval $n_s$ to the image at  $y^*=y^*(x',0,s^*)$ if $(n_s-1)t_1+t \leq T(y^*,\tilde{x}) \leq n_st_1+t -s^*$. Subject to this condition, the coordinate transform for image point $y^*$ is given by $\Tvk^{*(n_s)}(y^*)=(x_0',T(y^*,\tilde{x})+s^*)$ and the propagator matrices are given by $W^{s^*}$. Application of the box algorithm with $\Tvk^{*(n_s)}$ and $W^{s^*}$ with $\mathcal{A}$ incorporated in the amplitude $\tilde a^*$ yields the partial image
\begin{multline}
\label{equ:rtcimbox}
\Delta^{*(n_s,t+(n_s-1)t_1)}(y)= 
\sum_{\nu,k}\tilde a^*(y,\nu)\sum_{r=1}^{R_{\nu,k}} \alpha_{\nu,k}^{*(r)}(y)\\
\times \sum_{(\xi',\tau)} e^{\ii \langle \Tvk^{*(n_s)}(y),(\xi',\tau)\rangle}\hat g(\xi',\tau)\hat\beta_{\nu,k}(\xi',\tau) \hat\chi_{\nu,k}(\xi',\tau)\hat\vartheta_{\nu,k}^{*(r)}(\xi',\tau).
\end{multline}
In the case of conjugate points ($n_s>1$), we proceed with the evaluation of $H$ for the subsequent half wave equation reverse-time continuation of the wave fields $w_{r,+}^{(n_s)}(y,t+(n_s-1)t_1)$ (cf. Algo. \ref{tab:rtc}, Part II). Consider continuation of $w_{r,+}^{(n_s)}(y,t+(n_s-1)t_1)$ to $w_{r,+}^{(n_s)}(y,t+(n_s-2)t_1)$ ($n_p=2$). In this process, we compute the bicharacteristics $(\tilde y^s(y_0,\eta_0),\tilde \eta^s(y_0,\eta_0))$ for $s\in(0,t_1]$. We can easily monitor the condition $s^*=t+(n_s-1)t_1-T(\tilde y^{s*},\tilde{x})$ under which we obtain a contribution to the image at $\tilde y^{s*}$. The associated coordinate transform is given by $\tilde T_{\nu,k}^{*(n_s)}(\tilde y^{s*})=y_0$, and the propagator matrices by $\tilde W^*=W(y_0,\eta_0,s^*)$. With these ingredients, application of the box algorithm yields the partial image $\tilde \Delta^{*(n_s,n_p)}(y)$; we obtain the final image
$$
\Delta_{d_\Sigma}(y)=\sum_{n_s=1}^{N_s} \Delta^{*(n_s,t+(n_s-1)t_1)}(y) + \sum_{n_p=2}^{N_s}\sum_{n_s=n_p}^{N_s}\tilde \Delta^{*(n_s,n_p)}(y).
$$
The structure of the inverse scattering procedure is summarized in Algo. \ref{tab:rtcim}. Note that in the evaluation of the partial images $\Delta^*$ and $\tilde\Delta^*$, we can gather the 
incident angles $\eta(y^*;\nu,k;n_s,n_p)$ of the reverse-time continued wave field, which we can, for instance, use for monitoring scattering angles as illustrated in Section \ref{sec:angle}.

\begin{figure}[tb]
\centering
\includegraphics[width=\linewidth]{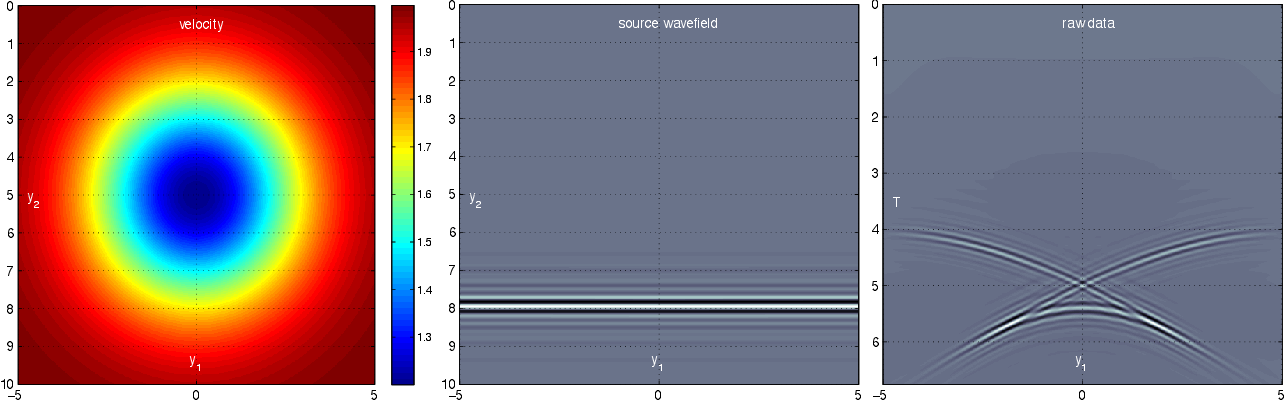}
\caption{\label{FigBP1} Reverse-time continuation from a boundary: Velocity model (left), initial wave field at $t=0$ (center) and data collected at the boundary $y_2=0$ (right).}\vspace{-5mm}
\end{figure}

\begin{figure}[htb]
\centering
\includegraphics[width=\linewidth]{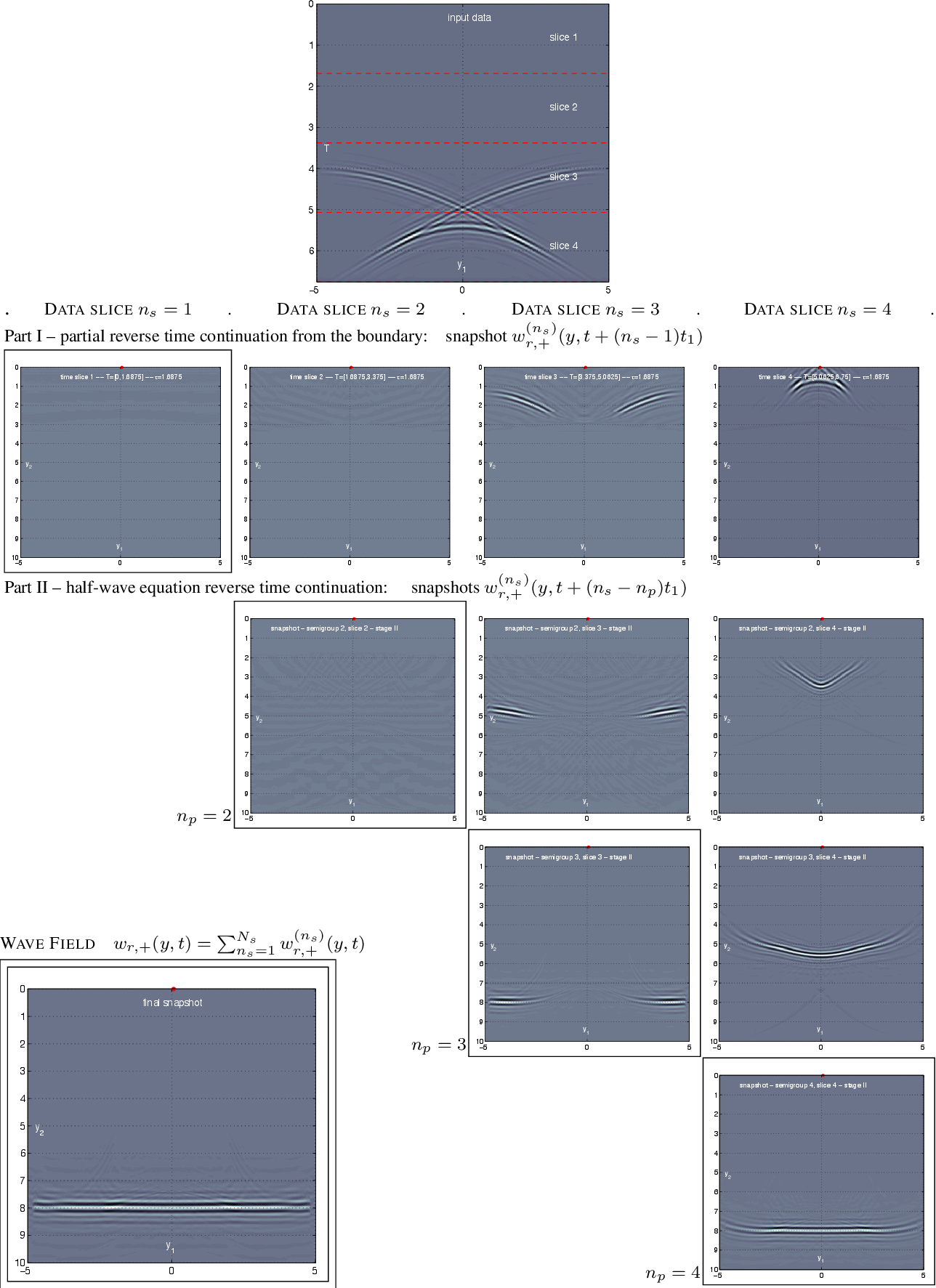}
\caption{\label{FigBPalgo} Top: Partitioning of data in Fig. \ref{FigBP1} (right) into $4$ time intervals. Center: Reverse-time continuation from the boundary, yielding $w_{r,+}^{(n_s)}(y,t+(n_s-1)t_1),\;n_s=1,\dots,4$ (top row); 
reverse-time continued wave fields $w_{r,+}^{(n_s)}(y,t+(n_s-n_p)t_1),\;n_s=n_p,\dots,4$ for $n_p=2,3$ and $4$, respectively (rows $2$ to $4$). Full reverse-time continued wave field $w_{r,+}(y,t=0)$ (bottom left corner).}
\end{figure}

\section{Numerical examples}
\label{sec:5}

We illustrate the performance of our algorithm in two applications:
Reverse-time continuation from the boundary of an upgoing wave field
in the presence of conjugate points and imaging of conormal
singularities using reverse-time continuation of boundary reflection
data. We consider the isotropic case. Although applicable in general
dimension, we restrict ourselves here to dimension $2$. The sources in
these examples, $g(x',t)$ and $d_\Sigma(x',t)$, respectively, are
generated using a time domain finite difference method. The
computational domain is of size $N \times N = 512 \times 512$.

\subsection{Reverse-time continuation from a boundary in the presence of caustics}

Here, we illustrate reverse-time continuation of boundary data
generated by a horizontal plane wave traveling upward through a low
velocity lens. The model is plotted in Fig. \ref{FigBP1} (left) and
consists of a Gaussian low wave speed lens with a variation of
40\% of the peak wave speed of the background model. The initial wave
field at $t=0$ is plotted in Fig. \ref{FigBP1} (center) and the
generated boundary data at $y_2=0$ are plotted in Fig. \ref{FigBP1}
(right).

In Fig. \ref{FigBPalgo} (top), we plot the data obtained after de-recomposition of the time domain finite difference data in Fig. \ref{FigBP1} (right) using the wave packet transform. Note that in this de-recomposition step, we can perform denoising, data regularization, or phase-space filtering (dip angle, wave number, location) in the wave packet domain and initiate ``beams'' \cite{BrandsbergEtgen2003}. 
We set $t=0$ and 
monitoring of the points of rank-deficiency of the matrix $W_1^t$ reveals that we
need $N_s=4$ time intervals and hence a total of three semigroup decompositions to avoid the formation of caustics in each step of the partial reverse-time continuation. 
The partitioning of the data in four time slices is indicated with red dashed lines in Fig. \ref{FigBPalgo} (top).

\clearpage

The center plots in Fig. \ref{FigBPalgo} show the partial outputs of
the reverse-time continuation procedure described in Section
\ref{sec:comprtc} and illustrate its logic and structure. Each column
corresponds with one time interval of the data (from left to right,
data slice $n_s=1,\dots,4$, respectively), and transition from row $i$
to row $i+1$ corresponds with a semigroup re-decomposition and
subsequent half wave equation reverse-time continuation step: The top
row plots the wave fields
$w_{r,+}^{(n_s)}(y,t+(n_s-1)t_1),\;n_s=1,\dots,4$, and the second, third,
and last rows show
$w_{r,+}^{(n_s)}(y,t+(n_s-n_p)t_1),\;n_s=n_p,\dots,4$ for $n_p=2,3$
and $4$, respectively (the reverse-time continued wave fields
$w_{r,+}^{(n_s)}(y,t=0)$ obtained for the four time intervals are
indicated by black frames).  The final reverse-time continued wave
field $w_{r,+}(y,t=0)=\sum_{n_s=1}^{N_s}w_{r,+}^{(n_s)}(y,t=0)$ is
plotted in the bottom left corner of Fig. \ref{FigBPalgo} (black solid
double-frame) and reproduces well the initial wave field at time $t=0$
(cf. Fig. \ref{FigBP1} (center)). Despite several discrete wave packet
transform re-decomposition steps involved in computing the reverse-time continuation (semigroup), the amplitude is accurate. In
particular, we note that the edges of the cusp in the data are well
focussed.

Note that time intervals $1$ and $2$ do not contain any significant
energy. With the proposed procedure, it is possible to compute only
the wave field for time intervals $3$ and $4$ (requiring no
computation time and memory for time slices 1 and 2). Time intervals
$1$ and $2$ have nonetheless been included in Fig. \ref{FigBPalgo} for
completeness of the presentation.

\begin{figure}[tb]
\centering
\includegraphics[width=\linewidth]{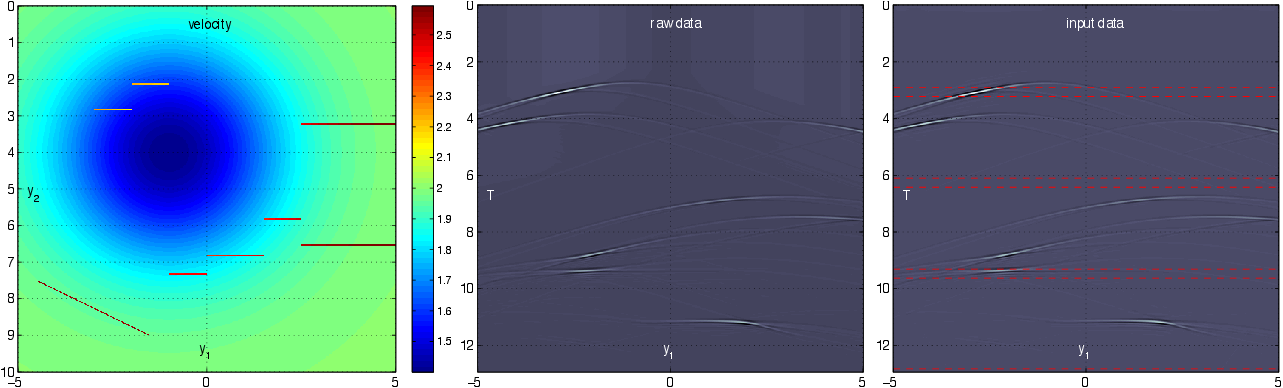}
\caption{\label{FigBPIMmodel} Imaging of conormal singularities: velocity model with line reflectors (left); time domain finite difference snapshot for source position $\tilde x=0$ (center); data after wave packet de-recomposition (right), the overlap of the time intervals partitioning the data are indicated by red dashed lines.}
\end{figure}

\subsection{Imaging of conormal singularities}

\begin{figure}[tb]
\centering
\includegraphics[width=\linewidth]{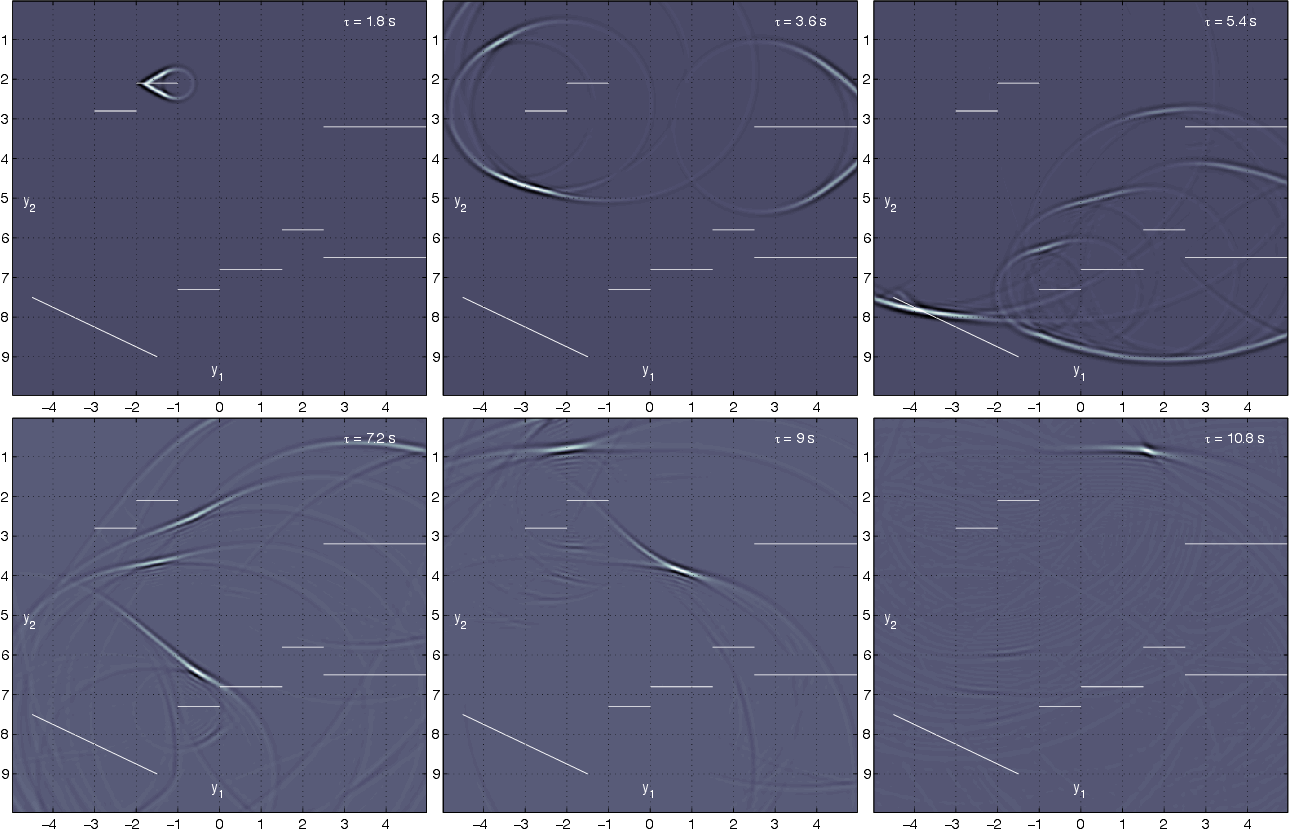}
\caption{\label{FigBPIMsnaps} The wave field generated by a source at $\tilde x=0$ for different time instances after subtraction of the reference wave field which is obtained when the line reflectors (indicated in white) are not present.} \vspace{-2mm}
\end{figure}

\begin{figure}[tb]
\centering
\includegraphics[width=0.89\linewidth]{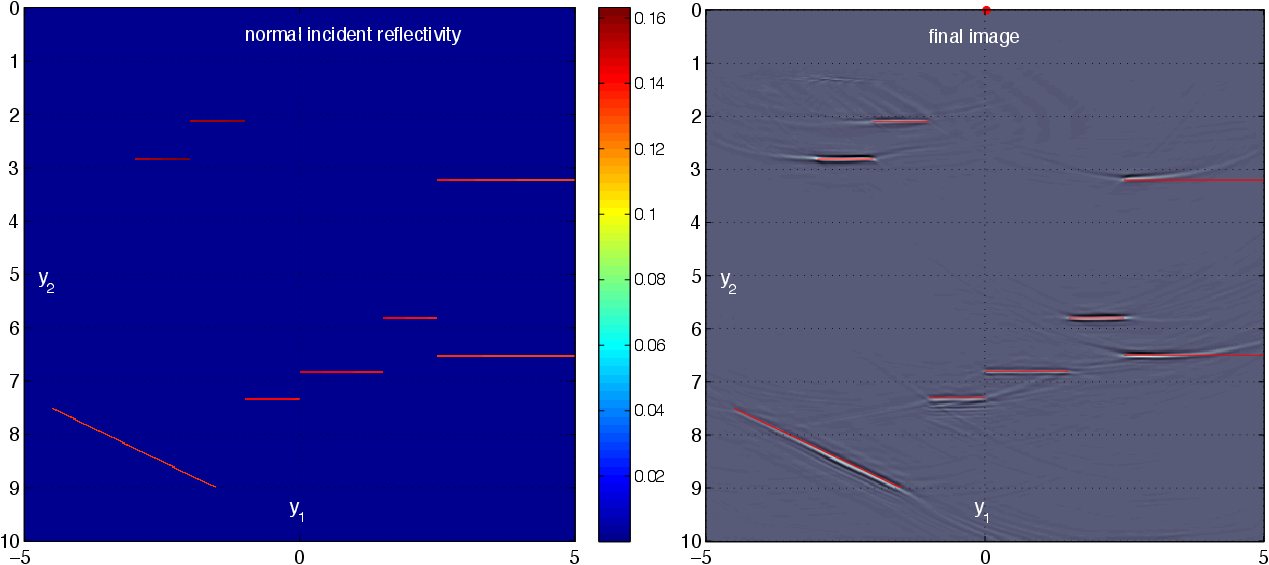}\vspace{-3mm}
\caption{\label{FigBPIMimage} Normal incident reflectivity of the model (left) and image $\Delta_{d_\Sigma}(y)$ (right).}\vspace{-2mm}
\end{figure}

We proceed with a numerical illustration of imaging of conormal
singularities by reverse-time continuation from the boundary using the
wave packet based computational procedure developed in Section
\ref{sec:4}. The velocity model is plotted in Fig. \ref{FigBPIMmodel}
(top left). It consists of a decentered Gaussian low velocity (30\%
peak contrast with respect to the background velocity) and contains
several horizontal line reflectors and one deep tilted line
reflector. The (normal incident) reflectivity of the line reflectors
varies with location and is documented in Fig. \ref{FigBPIMimage}
(left). The data are generated using time domain finite difference and
a Ricker wavelet with a peak frequency of $7Hz$. The single source is
located at the center of the boundary, $\tilde x=0$. In
Fig. \ref{FigBPIMsnaps}, we plot the wave field generated in the
subsurface for several time instances (for better visibility, we
substracted the wave field that is obtained when the line reflectors
are not present). Despite the simplicity of the model, we observe a
relatively complicated wave field and, for late time instances, the
formation of caustics. Also note that artifacts from nonperfectly
absorbing boundaries and from multiple reflections, and in particular
some numerical dispersion at large times are present in the simulated
wave field and consequently also in the data, which we plot in
Fig. \ref{FigBPIMmodel} (center).  The data $d_\Sigma(x',t)$ after
de-recomposition using the discrete wave packet transform are plotted
in Fig. \ref{FigBPIMmodel} (right). During this de-recomposition step,
we can also regularize and preprocess the data (denoising, phase
space filtering).


In this example, we need $N_s=4$ time intervals to avoid conjugate
points within each propagation step in the computational procedure
described in Section \ref{sec:4} and outlined in
Algo. \ref{tab:rtcim}. This partitioning into time intervals 
is detected numerically from the points of rank-deficiency of the matrix $W_1^t$ of the Hamiltonian system as detailed in \cite{dHUVW:2013} and
indicated in Fig. \ref{FigBPIMmodel} (right).

\begin{figure}[htb]
\centering
\includegraphics[width=\linewidth]{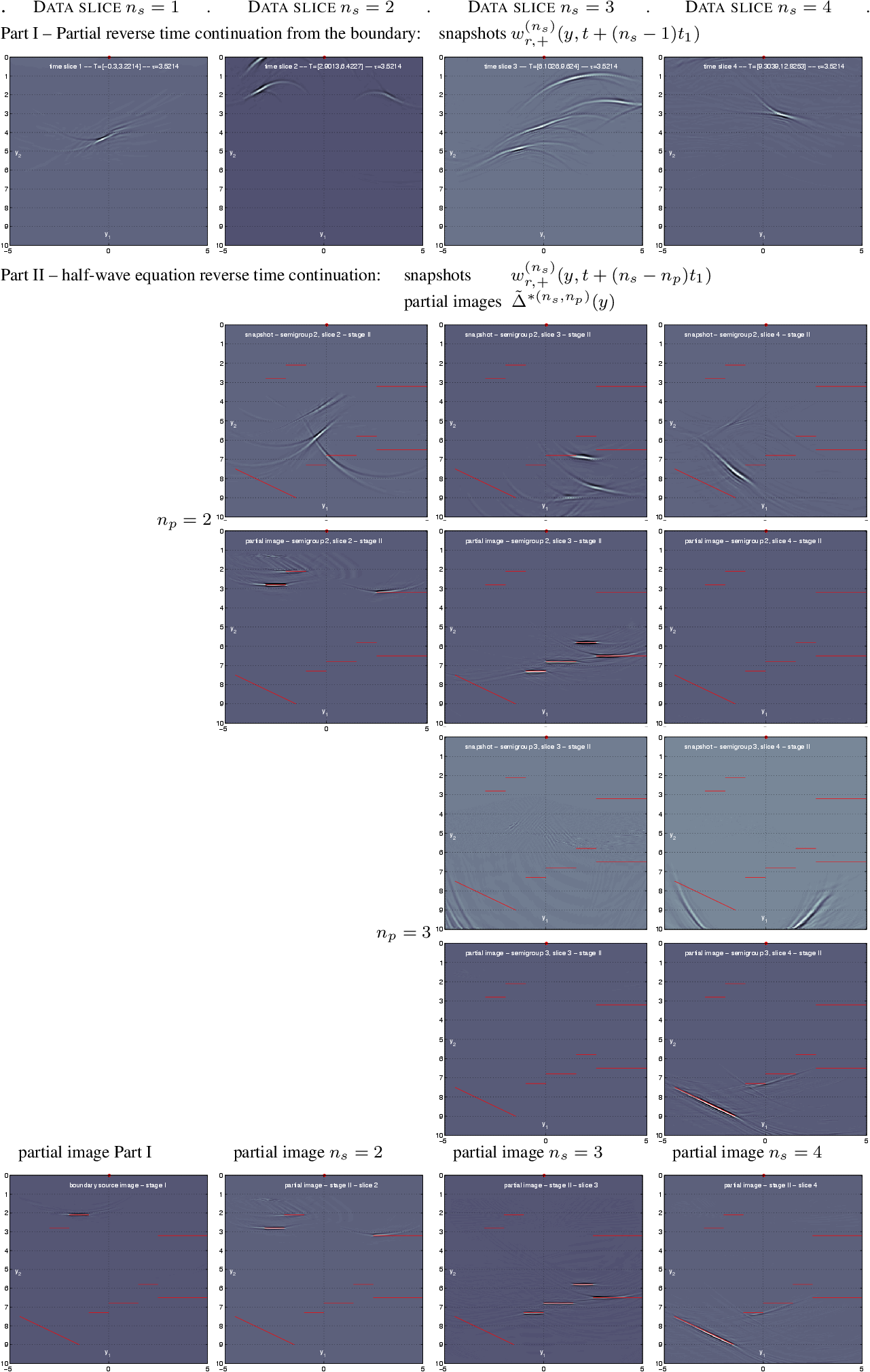}
\caption{\label{FigBPIMalgo} Partial reverse-time continuation from the boundary of the four time intervals in Fig. \ref{FigBPIMmodel} (right); reverse-time continuation and imaging following a semigroup re-decomposition of $w_{r,+}^{(n_s)}(y,t+(n_s-n_p+1 )t_1)$  for $n_p=2,\,n_s=2,\dots,4$ (center top rows) and for $n_p=2,\,n_s=2,\dots,4$ (center top rows):  snapshots $w_{r,+}^{(n_s)}(y,t+(n_s-n_p)t_1)$ and partial images $\tilde \Delta^{*(n_s,n_p)}(y)$; partial images produced by Part I (bottom row, left) and by Part II for $n_s=2,3,4$ (bottom right).}
\end{figure}

We approximate the source signature with a single delta function at its temporal maximum and compute the source wave field by evaluating the Hamiltonian and Hamilton-Jacobi equations (dynamic ray tracing). 
In Fig. \ref{FigBPIMalgo}, the partial images and reverse-time continued wave fields produced by the procedure described in Section \ref{sec:4} are plotted, organized according to its hierarchical semigroup structure (cf. Fig. \ref{FigBPalgo}). 
Each column corresponds with one time interval of the data (from left to right, data slice $n_s=1,\dots,4$, respectively), and transition from (group of) row(s) $i$ to (group of) row(s) $i+1$ corresponds with a semigroup re-decomposition and subsequent half wave equation reverse-time continuation and partial imaging step.
The top row shows the snapshots $w_{r,+}^{(n_s)}(y,t+(n_s-1)t_1),\;n_s=1,\dots,4$ produced by partial reverse-time continuation from the boundary of the $4$ data slices (Algo. \ref{tab:rtcim}, Part I). The corresponding partial image $\sum_{n_s}\Delta^{*(n_s,t+(n_s-1)t_1)}(y)$ obtained during this step is plotted in the bottom left corner of Fig. \ref{FigBPIMalgo}. At this stage, data slice 1 is fully reverse-time continued ($t=0$) while data slices 2 to 4 will be further reverse-time continued after a semigroup re-decomposition (and enter Part II of Algo. \ref{tab:rtcim}). The second and third groups of rows plot the output of Part II (cf. Algo. \ref{tab:rtcim}) of the procedure for $n_p=2$ and $n_p=3$, respectively: $w_{r,+}^{(n_s)}(y,t+(n_s-n_p)t_1)$ (top rows) and $\tilde \Delta^{*(n_s,n_p)}(y)$ (bottom rows). We stop the semigroup iteration at $n_p=3$ because the energy of the data wave fields $w_{r,+}^{(n_s)}(y,t+(n_s-3)t_1)$ has already passed the image region of interest, and further reverse-time continuation would not add any energy to the final image. The partial image contributions of data slices 2 to 4 are plotted in the bottom row of Fig. \ref{FigBPIMalgo} (second to fourth columns).

Let us finally turn our attention to the image
$\Delta_{d_\Sigma}(y)=\sum_{n_s=1}^{4} \Delta^{*(n_s,t+(n_s-1)t_1)}(y) + \sum_{n_p=2}^{3}\sum_{n_s=n_p}^{4}\tilde \Delta^{*(n_s,n_p)}(y)$,
which is plotted in Fig. \ref{FigBPIMimage} (right). We observe that all the reflectors are imaged correctly and well focused, regardless of their depth, dip angle and background velocity. Note that we could further focus the image by using the full source signature instead of a delta source approximation. Certain reflectors are partially outside of the zone of illumination (e.g. the two rightmost reflectors at depths $y_2=3.2$ and $y_2=6.5$) and hence produce smiling ``tails'' caused by the truncation of the wave field in the data (cp. Kirchhoff migration). Similarly, the corners of the line reflectors act as point diffractors and produce tails according to partial illumination and restricted geometry. Note that the ringing artifacts in the data components stemming from the two deepest reflectors are also present in the image -- the algorithm produces an image from the data, with its imperfections. This is also the case for the artifact at depth $y_2=1.3$ in the image, which results from an imperfectly removed direct arrival (cf. Fig. \ref{FigBPIMmodel} (right), $(y_1,t)=(3,2.8)$).


\subsection{Restricted angle transform}
\label{sec:angle}

\begin{figure}[tb]
\centering
\includegraphics[width=0.99\linewidth]{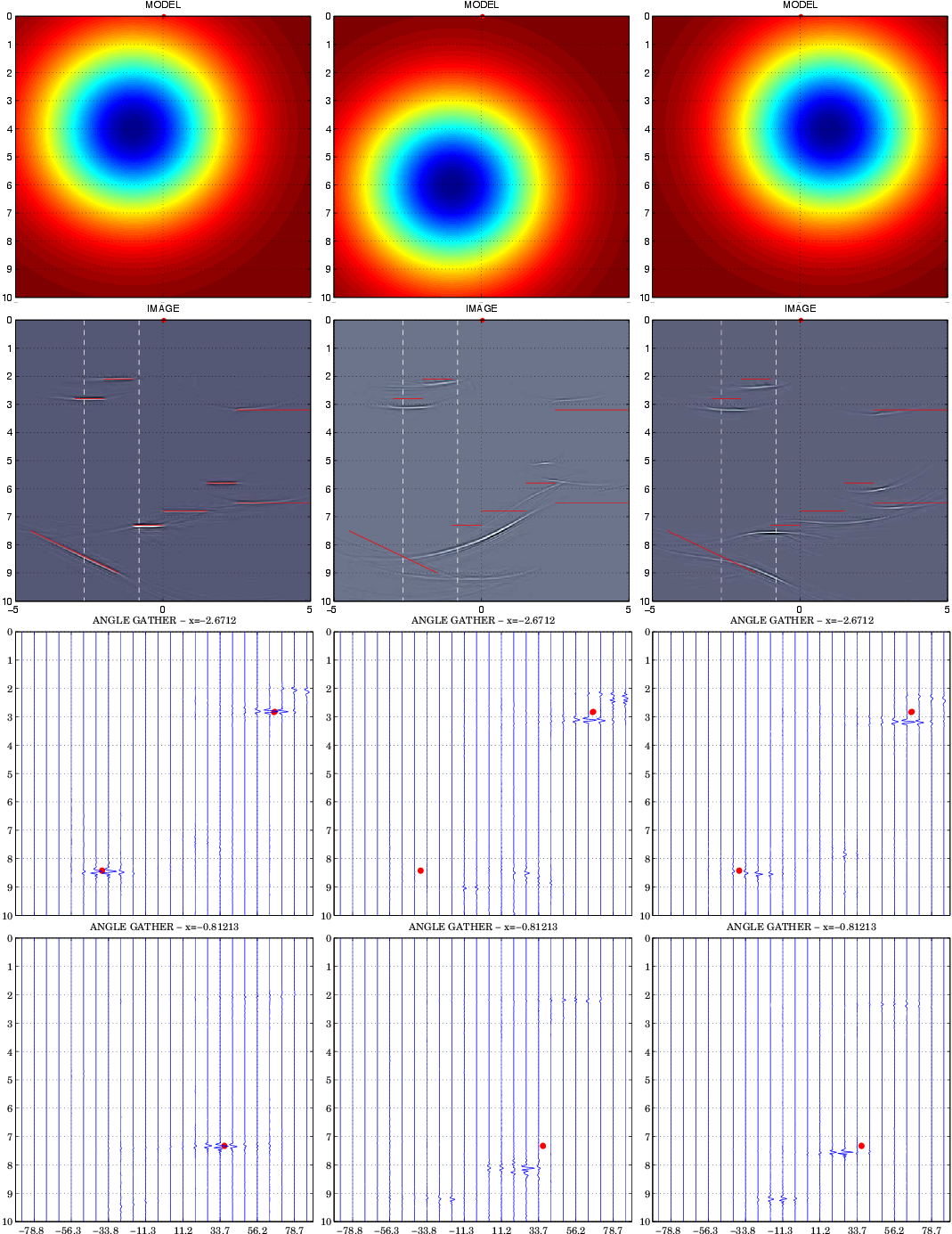}\vspace{-3mm}
\caption{\label{FigAngleGather} Velocity models (top row), resulting images (second row) and angle gathers at horizontal positions $x=-2.67$ (third row) and $x=-0.81$   (bottom row): correct velocity model (left column) and wrong velocity models (center and right column). The red dots indicate the specular reflection points for the true velocity model.\vspace{-4mm}}
\end{figure}

\begin{figure}[tb]
\centering
\includegraphics[width=0.99\linewidth]{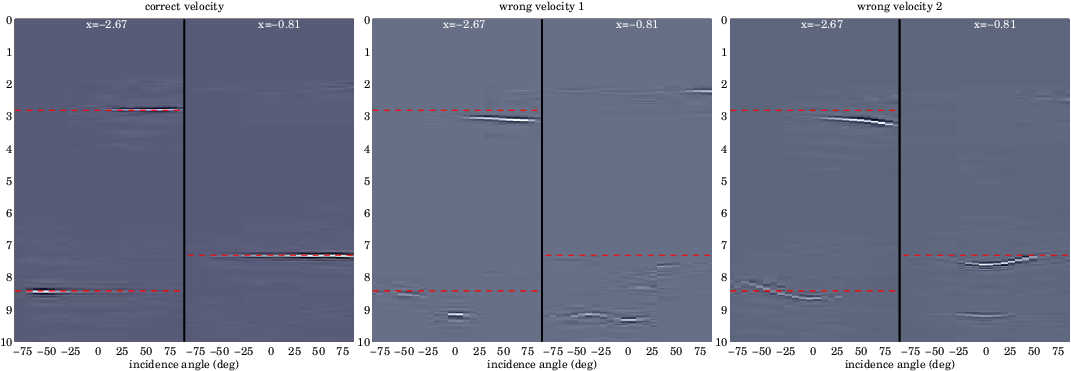}\vspace{-3mm}
\caption{\label{FigImGather} Image gathers for horizontal positions $x=-2.67$ and $x=-0.81$ as considered in Fig. \ref{FigAngleGather} (left and right half of images, respectively): correct velocity model (left) and wrong velocity models (center and right, respectively). The red dashed lines indicate the depth of the line reflectors.\vspace{-3mm}}
\end{figure}

Since the proposed boundary source reverse-time continuation and
imaging procedures rely on the dyadic parabolic decomposition, angular
information can be extracted for the source and scattered wave
fields. Indeed, for a given frequency box $\hat\chi_{\nu,k}$, the
incidence angles of the wave fronts are known in each image
point. This information can be directly incorporated into the imaging
process. Indeed, we can directly generate so-called image gathers in
incidence angles (which can be converted to scattering angles), that
is, generate images for particular incidence angles. This is
illustrated in Figs. \ref{FigAngleGather} and \ref{FigImGather}.

In Fig. \ref{FigAngleGather} (second row), we display the images obtained for a single source with the correct velocity model (left column), as in the previous section, and with two wrong velocity models (center and right column, respectively; the corresponding velocity models are plotted in the first row of Fig. \ref{FigAngleGather}).

\clearpage

In the third and fourth rows of Fig. \ref{FigAngleGather}, the images are analyzed as a function of incidence angle at the image points for two distinct boundary locations, respectively.
Geometrically, the image of a singularity at one surface location is significant at one incidence angle only; we indicate the incidence angle at specular reflection with a red dot.
The finite-frequency content of the wave packets results in a slight spread around these specular reflection angles. For the correct background model, the images are in phase at the depth of the reflector in the neighborhood of the specular reflection angles, while they are not if the wrong velocity model is used.

We evaluate images of the singularities for multiple sources
and rearrange them in terms of local incidence angle (image gather) \cite{Exxon,
  Yoon-Guo-Cai-Wang_2011, Yoon-Marfurt-Starr_2004,
  Zhang-McMechan_2011}. 
In Fig. \ref{FigImGather}, the images obtained using the correct (left column) and two wrong (second and third column) velocity models are plotted (the corresponding velocity models are plotted in the first row of Fig. \ref{FigAngleGather}). In case the correct background velocity
function is used, up to illumination effects, the images generated at
different angles are the same; this reflects a redundancy in such
data. If we perturb the background by moving the smooth lens, we still
obtain coherent images; however, the singularities move with
changing incidence angle. This behavior can be exploited
to develop a procedure for reflection tomography \cite{Stolk-deHoop,
  Burdick-etal_2013}.

\section{Discussion}
\label{sec:6}

We obtained a representation of RTM in terms of a FIO associated with
a canonical graph. We then developed a fast algorithm for solving the
wave equation with a boundary source and homogeneous initial
conditions using the dyadic parabolic decomposition of phase space,
adapting our algorithm for the computation of the action of FIOs associated with canonical graphs
\cite{AndersondHW:2012}, which is the key component of its
application. We admit the formation of caustics.

Our algorithm is organized by frequency boxes $\hat\chi_{\nu,k}$
following the dyadic parabolic decomposition of phase space and yields
accuracy $\mathcal{O}(2^{-k/2})$ at frequency scale $k$.  We obtain an
effective one-step multiscale procedure for reverse-time continuation
from the boundary for a given time interval, from $T_1$ to $t$,
say. In this process, we can apply the imaging condition and obtain a
reverse-time-migration imaging algorithm.

While numerical illustrations have been devised here for dimension
$2$, the concepts and computational procedures are valid for arbitrary
dimension.

In the presence of conjugate points, we split the time interval for reverse-time continuation into a sequence of smaller time intervals and reverse-time continue partial wave fields subsequently for these time intervals using the semigroup property of the RTM operator. Numerically, this implies one discrete wave packet transform re-decomposition of the wave fields for each transition point from one time interval to another. After the first semigroup re-decomposition, reverse-time continuation essentially reduces to the evaluation of the wave equation for the propagation of an initial wave field, and any of the algorithms developed in \cite{AndersondHW:2012} could be used as a computational basis. Here, we proposed a ``box algorithm'' due to its favorable computational complexity and practical accuracy.

The computational complexity of our algorithm is of the order $\mathcal{O}(N^n\log(N))$
per frequency box for each semigroup step for an $n$-dimensional grid of side length $N$. It arises essentially from
the complexity of the unequally spaced FFTs involved in the box
algorithm (cf. \cite{AndersondHW:2012} for details). Computations for
each individual box are independent and hence embarrassingly parallel.
The computational cost of RTM imaging is roughly
twice that of reverse-time continuation of the wave field from the
boundary because of the additional unequally spaced FFTs that need to
be evaluated for producing the image. Note that with the exception of
the source wave field travel times and amplitudes and one single
snapshot during each semigroup re-decomposition, our procedure does
not require the computation and storage of snapshots. Its
computational and memory requirements are therefore of the order of
the one-step evaluation of Cauchy initial value problems for evolution
equations in \cite{AndersondHW:2012}.

Evaluation of the RTM operator for all wave packets associated with a given frequency boxes $\hat\chi_{\nu,k}$ at once requires the existence of a homogeneous boundary layer near the acquisition surface. When the wave speed is not constant near the boundary, we need to localize computations and either introduce a partitioning of the acquisition surface or use wave packets as individual local data quanta, the latter yielding a wave packet based procedure at the price of increased computational complexity with respect to a frequency box driven algorithm.

The total number of frequency boxes is $\mathcal{O}(N^{(n-1)/2})$, this number can be slightly reduced by replacing frequency boxes (tiles) with wedges as in \cite{Candes2009}, yet at the price of losing the connection to the data wave packets. Depending on the data and the imaging target, not all boxes need to be computed. 
Indeed, our algorithm provides phase-space localized control for the data (scale, orientation, position of the data wave packets) as well as the image (scale, orientation and position of reverse-time continued data wave packets; full angular information such as scattering angle, and reflector dip angle). 
Together with the fact that only a few time steps need to be computed for producing an image, this makes our algorithm particularly attractive for (limited aperture) array data, partial imaging and target-oriented imaging. An additional asset of our approach is that incident angles of wave fronts can also be imaged, enabling the straightforward application of restricted angle transforms.
 
We note that by viewing wave packets as localized plane waves, our method can be related to plane-wave and beam-wave migration \cite{BrandsbergEtgen2003}. Here, we can construct ``beams'' as reverse-time continued data wave packets based on phase-space localized paraxial approximation in geodesic coordinates. 
In the context of imaging with ambient noise using body waves and beamforming \cite{poli2012body}, one generates a cross correlation matrix between two distant receiver arrays at which ambient noise generated data are obtained, and one takes inner products with wave packets and can enhance particular wave constituents prior to applying the imaging operator.

Reverse-time continuation from the boundary can in principle be
generalized to extended imaging using multisource data based on
downward continuation \cite{Stolk2005a}. The corresponding evolution
equation replacing (2.6) can be found in
\cite[Eq. (17)]{Duchkov2010}. In this case, the evolution equation is
defined in $(2n-1)$-dimensional extended space.


\end{document}